# LARGE DEVIATIONS AND LAWS OF THE ITERATED LOGARITHM FOR THE LOCAL TIMES OF ADDITIVE STABLE PROCESSES[1]

By Xia Chen

*University of Tennessee*

We study the upper tail behaviors of the local times of the additive stable processes. Let $X_1(t), \ldots, X_p(t)$ be independent, $d$-dimensional symmetric stable processes with stable index $0 < \alpha \leq 2$ and consider the additive stable process $\overline{X}(t_1, \ldots, t_p) = X_1(t_1) + \cdots + X_p(t_p)$. Under the condition $d < \alpha p$, we obtain a precise form of the large deviation principle for the local time

$$\eta^x([0,t]^p) = \int_0^t \cdots \int_0^t \delta_x(X_1(s_1) + \cdots + X_p(s_p)) \, ds_1 \cdots ds_p$$

of the multiparameter process $\overline{X}(t_1, \ldots, t_p)$, and for its supremum norm $\sup_{x \in \mathbb{R}^d} \eta^x([0,t]^p)$. Our results apply to the law of the iterated logarithm and our approach is based on Fourier analysis, moment computation and time exponentiation.

**1. Introduction.** Throughout, $X_1(t), \ldots, X_p(t)$ are independent $d$-dimensional symmetric stable processes with identical distribution. We use the notation $X(t)$ for a stable process with the same distribution as $X_1(t), \ldots, X_p(t)$. In this paper, the stable index $\alpha \in (0, 2]$. By our assumptions, there is a continuous function $\psi(\lambda) \geq 0$ on $\mathbb{R}^d$ with

$$\psi(r\lambda) = r^\alpha \psi(\lambda) \quad \text{and} \quad \psi(-\lambda) = \psi(\lambda), \qquad r > 0, \ \lambda \in \mathbb{R}^d,$$

such that

$$(1.1) \qquad \mathbb{E} e^{i\lambda \cdot X(t)} = e^{-t\psi(\lambda)}, \qquad t \geq 0, \ \lambda \in \mathbb{R}^d.$$

Since we only consider nondegenerate stable processes, there is a constant $C > 0$ such that

$$C^{-1}|\lambda|^\alpha \leq \psi(\lambda) \leq C|\lambda|^\alpha.$$

---

Received May 2005; revised March 2006.

[1]Supported in part by NSF Grant DMS-04-05188.

Unless assuming otherwise, $X_1(0) = \cdots = X_p(0) = 0$.

The following $p$-parameter, $d$-dimensional random field:

$$\overline{X}(t_1, \ldots, t_p) = X_1(t_1) + \cdots + X_p(t_p), \qquad (t_1, \ldots, t_p) \in (\mathbb{R}^+)^p,$$

is called an additive stable process.

Since they locally resemble stable sheets, and since they are more amenable to analysis, additive stable processes first arose to simplify the study of stable sheets (see [9, 10, 17] and [18]). They also arise in the theory of intersection and self-intersection of stable processes (see [15, 22] and [26]). As pointed out below [see (1.12)], the local time of additive processes is actually an intersection local time as $p = 2$. We refer the reader to [1, 2, 4, 6, 7, 8] and [25] for some recent developments in the large deviations for intersection local times. We also point out the reference [5] for the study on the small ball probabilities of the additive stable processes. The study of additive processes also connects to probabilistic potential theory. We mention [16, 20, 21, 22] and refer the reader to the detailed discussion and for further reference.

In this work, we consider the local times of $\overline{X}(t_1, \ldots, t_p)$ which are formally given as

$$\eta^x(I) = \int_I \delta_x(X_1(s_1) + \cdots + X_p(s_p)) \, ds_1 \cdots ds_p, \qquad x \in \mathbb{R}^d, \ I \subset (\mathbb{R}^+)^p.$$

We rely on two recent papers by Khoshnevisan, Xiao and Zhong [23, 24] for the constructions of the local time $\eta^x(I)$. In their papers, Khoshnevisan, Xiao and Zhong [23, 24] consider a more general multiparameter random field named additive Lévy process, which is generated by independent Lévy processes. In their construction, $\eta^x(I)$ is defined as the density function of the occupation measure $\mu_I$:

$$\mu_I(A) = \int_I \delta_{\overline{X}(s_1, \ldots, s_p)}(A) \, ds_1 \cdots ds_p, \qquad A \subset \mathbb{R}^d,$$

in the case when $\mu_I$ is absolutely continuous with respect to the Lebesgue measure on $\mathbb{R}^d$. Applying Theorem 1.1 in [23] to our setting, the local time $\eta^x(I)$ exists for every super interval $I \subset (\mathbb{R}^+)^p$ if and only if

$$(1.2) \qquad\qquad\qquad\qquad d < \alpha p.$$

Under (1.2),

$$(1.3) \qquad\qquad\qquad \int_{\mathbb{R}^d} [\eta^x(I)]^2 \, dx < \infty \qquad \text{a.s.}$$

for every finite $d$-dimensional interval $I \subset (\mathbb{R}^+)^p$ (Theorem 1.3 of [23]). Further, (1.2) also implies that almost surely, the local time

$$\eta^x([0, t]^p), \qquad (x, t) \in \mathbb{R}^d \times \mathbb{R}^+,$$

is jointly continuous in $(x, t)$ (Corollary 3.3 of [24]).

We mention that in the stable case, Khoshnevisan, Xiao and Zhong ([24], Theorems 4.3 and 5.3) carried out some tail estimates for the local time $\eta$ which yields a sharp rate.

In this paper, our goal is to establish the large deviations and the laws of the iterated logarithm for the local times of additive stable processes. In particular, we shall identify, as much as we can, the constants appearing in these limit forms.

Recall that the characteristic exponent $\psi(\lambda)$ is defined by (1.1) and write

$$(1.4) \qquad \rho = \sup_{\|f\|_2 = 1} \int_{\mathbb{R}^d} \left[ \int_{\mathbb{R}^d} \frac{f(\lambda + \gamma) f(\gamma)}{\sqrt{1 + \psi(\lambda + \gamma)} \sqrt{1 + \psi(\gamma)}} \, d\gamma \right]^p d\lambda,$$

where

$$\|f\|_2 = \left( \int_{\mathbb{R}^d} f^2(\lambda) \, d\lambda \right)^{1/2}.$$

Clearly, $\rho > 0$. We now prove that under the condition (1.2), $\rho < \infty$. Indeed, by Hölder inequality

$$\left[ \int_{\mathbb{R}^d} \frac{f(\lambda + \gamma) f(\gamma)}{\sqrt{1 + \psi(\lambda + \gamma)} \sqrt{1 + \psi(\gamma)}} \, d\gamma \right]^p$$
$$\leq \left( \int_{\mathbb{R}^d} |f(\lambda + \gamma) f(\gamma)| \, d\gamma \right)^{p-1} \int_{\mathbb{R}^d} \frac{|f(\lambda + \gamma) f(\gamma)|}{(1 + \psi(\lambda + \gamma))^{p/2} (1 + \psi(\gamma))^{p/2}} \, d\gamma.$$

By the Cauchy–Schwarz inequality and shift-invariance,

$$\int_{\mathbb{R}^d} |f(\lambda + \gamma) f(\gamma)| \, d\gamma \leq \int_{\mathbb{R}^d} f^2(\gamma) \, d\gamma = 1.$$

Hence,

$$\int_{\mathbb{R}^d} \left[ \int_{\mathbb{R}^d} \frac{f(\lambda + \gamma) f(\gamma)}{\sqrt{1 + \psi(\lambda + \gamma)} \sqrt{1 + \psi(\gamma)}} \, d\gamma \right]^p d\lambda$$
$$\leq \int_{\mathbb{R}^d} \left( \int_{\mathbb{R}^d} \frac{|f(\lambda + \gamma) f(\gamma)|}{(1 + \psi(\lambda + \gamma))^{p/2} (1 + \psi(\gamma))^{p/2}} \, d\gamma \right) d\lambda$$
$$= \left[ \int_{\mathbb{R}^d} \frac{|f(\lambda)|}{(1 + \psi(\lambda))^{p/2}} \, d\lambda \right]^2 \leq \int_{\mathbb{R}^d} \frac{1}{(1 + \psi(\lambda))^p} \, d\lambda.$$

Thus,

$$(1.5) \qquad \rho \leq \int_{\mathbb{R}^d} \frac{1}{(1 + \psi(\lambda))^p} \, d\lambda < \infty$$

where the last step follows from (1.2).

Our first main theorem is the large deviation principle for $\eta^0([0,t]^p)$. By the scaling property of the stable processes $X_1(\cdot), \ldots, X_p(\cdot)$, it can be verified that

$$(1.6) \qquad \eta^0([0,t]^p) \overset{d}{=} t^{(\alpha p - d)/\alpha} \eta^0([0,1]^p).$$

Without loss of generality, we need only to consider $\eta^0([0,1]^p)$ instead in the following theorem.

THEOREM 1.1. *Under* (1.2),

$$(1.7) \quad \lim_{t \to \infty} t^{-\alpha/d} \log \mathbb{P}\{\eta^0([0,1]^p) \geq t\} = -(2\pi)^\alpha \frac{d}{\alpha} \left(1 - \frac{d}{\alpha p}\right)^{(\alpha p - d)/d} \rho^{-\alpha/d}$$

*where* $\rho$ *is given in* (1.4).

We now connect Theorem 1.1 with some known results. As $p = 1$ and $\alpha > 1$, we have that $\psi(\lambda) = c|\lambda|^\alpha$ and that

$$\rho = \sup_{\|f\|_2 = 1} \left[ \int_{-\infty}^{\infty} \frac{f(\lambda)}{\sqrt{1 + c|\lambda|^\alpha}} \, d\lambda \right]^2 = \int_{-\infty}^{\infty} \frac{1}{1 + c|\lambda|^\alpha} \, d\lambda.$$

Theorem 1.1 becomes a classic large deviation result for the local time of the stable process $X(t)$ (see, e.g., [27]):

$$\lim_{t \to \infty} t^{-\alpha} \log \mathbb{P}\{\eta^0([0,1]) \geq t\} = -(2\pi)^\alpha \frac{1}{\alpha} \left(1 - \frac{1}{\alpha}\right)^{\alpha - 1} \rho^{-\alpha}.$$

We mention the large deviations for the intersection local time formally given as

$$\alpha([0,t]^p) = \int_{\mathbb{R}^d} \left[ \prod_{j=1}^{p} \int_0^t \delta_x(X_j(s)) \, ds \right] dx$$

under the condition

$$(1.8) \qquad p(d - \alpha) < d.$$

We refer the reader to the recent papers [4, 6] and [7] for the details on this subject. In particular, as $p = 2$, in which case (1.2) and (1.8) are equivalent to "$d < 2\alpha$," we have (Theorem 1 in [7]) that

$$(1.9) \qquad \lim_{t \to \infty} t^{-\alpha/d} \log \mathbb{P}\{\alpha([0,1]^2) \geq t\} = -\frac{d}{\alpha} \left(\frac{2 - \alpha}{2\alpha M_\psi}\right)^{(2\alpha - d)/d},$$

where

$$(1.10) \quad M_\psi = \sup_{g \in \mathcal{F}_\psi} \left\{ \left( \int_{\mathbb{R}^d} |g(x)|^4 dx \right)^{1/2} - \int_{\mathbb{R}^d} \psi(\lambda) |\widehat{g}(\lambda)|^2 \, d\lambda \right\},$$

$$(1.11) \quad \mathcal{F}_\psi = \left\{ g \in \mathcal{L}^2(\mathbb{R}^d); \; \|g\|_2 = 1 \text{ and } \int_{\mathbb{R}^d} \psi(\lambda) |\widehat{g}(\lambda)|^2 \, d\lambda < \infty \right\},$$

$$\widehat{g}(\lambda) = \int_{\mathbb{R}^d} g(x) e^{i\lambda \cdot x} \, dx, \qquad \lambda \in \mathbb{R}^d.$$

On the other hand, by the fact that $p = 2$ and $d < 2\alpha$,

$$(1.12) \qquad\qquad \alpha([0,1]^2) \stackrel{d}{=} \eta^0([0,1]^2).$$

Consequently, (1.9) is a direct corollary of Theorem 1.1, in the light of Lemma A.2 given in the Appendix.

By the continuity of $\eta^x([0,t]^p)$ when viewed as a function of $x$, and by the fact that $\eta^x([0,t]^p)$ is locally supported, we have that almost surely $\sup_{x \in \mathbb{R}^d} \eta^x([0,t]^p) < \infty$. By scaling, we have that for any $t > 0$,

$$(1.13) \qquad\qquad \sup_{x \in \mathbb{R}^d} \eta^x([0,t]^p) \stackrel{d}{=} t^{(\alpha p - d)/\alpha} \sup_{x \in \mathbb{R}^d} \eta^x([0,1]^p).$$

It has been known (see, e.g., [13] and [19]) that as $p = 1$, $\sup_{x \in \mathbb{R}^d} \eta^x([0,1]^p)$ has a tail behavior same as $\eta^0([0,1]^p)$. The following theorem claims that it remains true for $p \geq 2$.

THEOREM 1.2. *Under* (1.2),

$$(1.14) \quad \begin{aligned} &\lim_{t \to \infty} t^{-\alpha/d} \log \mathbb{P} \left\{ \sup_{x \in \mathbb{R}^d} \eta^x([0,1]^p) \geq t \right\} \\ &\qquad = -(2\pi)^\alpha \frac{d}{\alpha} \left( 1 - \frac{d}{\alpha p} \right)^{(\alpha p - d)/d} \rho^{-\alpha/d}, \end{aligned}$$

*where the constant $\rho$ is given in* (1.4).

Theorems 1.1 and 1.2 apply to the following law of the iterated logarithm.

THEOREM 1.3. *Assume* (1.2). *Then for any $x \in \mathbb{R}^d$,*

$$(1.15) \quad \begin{aligned} &\limsup_{t \to \infty} t^{-(\alpha p - d)/\alpha} (\log \log t)^{-d/\alpha} \eta^x([0,t]^p) \\ &\qquad = (2\pi)^{-d} \left( \frac{\alpha}{d} \right)^{d/\alpha} \left( 1 - \frac{d}{\alpha p} \right)^{-(p-d/\alpha)} \rho \qquad a.s. \end{aligned}$$

*and*

(1.16)

$$\limsup_{t \to \infty} t^{-(\alpha p - d)/\alpha} (\log \log t)^{-d/\alpha} \sup_{x \in \mathbb{R}^d} \eta^x([0, t]^p)$$

$$= (2\pi)^{-d} \left( \frac{\alpha}{d} \right)^{d/\alpha} \left( 1 - \frac{d}{\alpha p} \right)^{-(p - d/\alpha)} \rho \qquad a.s.$$

Our approach consists of three tools: time exponentiation, Fourier transformation and moment estimation. To outline some key ideas used in this paper, we first cite a lemma given in [25].

LEMMA 1.4 (Lemma 2.3 in [25]). *Let $Y$ be any nonnegative random variable and let $\theta > 0$ be fixed. Assume that*

(1.17)

$$\lim_{n \to \infty} \frac{1}{n} \log \frac{1}{(n!)^\theta} \mathbb{E} Y^n = -\kappa$$

*for some $\kappa \in \mathbb{R}$. Then we have*

(1.18)

$$\lim_{t \to \infty} t^{-1/\theta} \log \mathbb{P}\{Y \geq t\} = -\theta e^{\kappa/\theta}.$$

In their original statement, König and Mörters assume that $\theta > 0$ is an integer. By examining their proof, we find that $\theta$ can be any positive number.

Lemma 1.4 simply says that in order to have a tail estimate for a nonnegative random variable with certain precision, one needs to understand its high moment asymptotics. In Section 2 we first introduce a theorem (Theorem 2.1 below) without proof (which will be given in later sections) in which the high moment asymptotics are evaluated for the local time of additive stable process stopped at $p$ independent exponential times. Then we prove Theorem 1.1 based on Theorem 2.1. Although the scheme of time exponentiation has become standard in the area of limit theory since the remarkable work done by Darling and Kac [11], it is not usual to see such an idea being used in the context of multiparameter processes, at least not at the level of precision carried out in this work.

In Section 3 we prove the lower bound for Theorem 2.1. By Fourier transformation the moment of the local time (run up to exponential times) can be represented as an $\mathcal{L}_p$-norm. Then the lower bound follows from a simple argument via spectral theory.

The upper bound of Theorem 2.1 is much harder than the lower bound and needs a completely different treatment. In Section 4 we shall establish a discrete version (Theorem 4.1) of Theorem 2.1. The argument is combinatorial and is partially inspired by the pioneer work of König and Mörters [25] despite some essential differences between the situations faced by them and by us. We shall adopt a probabilistic approach to handle the moment asymptotics which is no longer a probabilistic problem. In Section 5 we complete the proof of the upper bound for Theorem 2.1. In this section we follow an interesting procedure of discretization by Fourier transform.

We prove Theorem 1.2 in Section 6 and Theorem 1.3 in Section 7. The proof relies on the exponential integrability of the local time (Lemma 6.1) under the Hölder norm and on some results established in the previous sections.

In the Appendix, we prove two analytic lemmas.

The central part of this work is Theorem 4.1 which is similar in spirit to Proposition 2.2 in [25] where the high moments of intersection local times are estimated. Here we compare the present paper with the one by König and Mörters [25]. A key ingredient in both works is to write the moments in terms of $L_p$-norms. In the case of intersection local times (studied by König and Mörters), the $L_p$-norm is related to Green's function; while in the case of the local times of additive processes, the $L_p$-norm is related to the Fourier transform of Green's function. In Proposition 2.2 in [25], the domain of intersection is limited to a compact set; while in our case the independent stable processes are allowed to interact at everywhere in $\mathbb{R}^d$. Consequently, compactification of the state space is one of several key issues addressed in our argument. Finally, both Proposition 2.2 in [25] and Theorem 4.1 are proved by combinatorial approaches and therefore both treatments contains a certain procedure of discretization. As to be pointed out at the beginning of Section 5 below, the classic procedure adapted in [25] is no longer working in our setting. Our way of discretization is based on some delicate properties of Fourier transformation.

We end this section with the following comment: the moment asymptotics linked to the weak convergence have been investigated extensively. We refer the interested reader to the survey paper by Fitzsimmons and Pitman [14] for an overview. In the study of the weak convergence, the power of the moment is often fixed. However, much less has been explored on the high moment asymptotics (where the power tends to infinity) which are usually linked to the large deviations through some general large deviation principles like Lemma 1.4. The study of high moment asymptotics has great potential in solving some hard problems on the large deviations, such as the large deviations for the intersection local times of general Markovian and Gaussian processes, and for the local times of some other multiparameter processes like stable sheets. It is too early to see a full scale of applications possibly brought by the research of high moment asymptotics; we leave it to future study.

**2. Time exponentiation.** In the rest of the paper, we introduce the notations $\tau_1, \ldots, \tau_p$ for independent exponential times with parameter 1, and $\Sigma_n$ for the set of all permutations on $\{1, \ldots, n\}$. We assume the independence between $\{\tau_1, \ldots, \tau_p\}$ and $\{X_1(t), \ldots, X_p(t)\}$. At first, we try to represent the $n$th moment of the random variable

$$\eta^0([0, \tau_1] \times \cdots \times [0, \tau_p])$$

in a reasonably nice form.

By Fourier transform, for any $t_1, \ldots, t_p \geq 0$,

$$\eta^0([0, t_1] \times \cdots \times [0, t_p])$$
$$= \frac{1}{(2\pi)^d} \int_{\mathbb{R}^d} \left[ \int_{\mathbb{R}^d} \eta^x([0, t_1] \times \cdots \times [0, t_p]) e^{i\lambda \cdot x} \, dx \right] d\lambda$$
$$= \frac{1}{(2\pi)^d} \int_{\mathbb{R}^d} d\lambda \int_0^{t_1} \cdots \int_0^{t_p} \exp\{i\lambda \cdot (X_1(s_1) + \cdots + X_p(s_p))\} \, ds_1 \cdots ds_p$$

where the second step follows from the definition of the local times as the density of occupation measures. Hence, for any integer $n \geq 1$,

$$\mathbb{E}[\eta^0([0, t_1] \times \cdots \times [0, t_p])^n]$$
$$= \frac{1}{(2\pi)^{dn}} \int_{(\mathbb{R}^d)^n} d\lambda_1 \cdots d\lambda_n$$
$$\times \prod_{j=1}^p \int_{[0, t_j]^n} \mathbb{E} \exp\left\{ i \sum_{k=1}^n \lambda_k \cdot X(s_k) \right\} ds_1 \cdots ds_n.$$

Let $\sum_n$ be the permutation group on the set $\{1, \ldots, n\}$. By time rearrangement and by independence of the increments,

$$\int_{[0, t_j]^n} \mathbb{E} \exp\left\{ i \sum_{k=1}^n \lambda_k \cdot X(s_k) \right\} ds_1 \cdots ds_n$$
$$= \sum_{\sigma \in \Sigma_n} \int_{\{0 \leq s_1 \leq \cdots \leq s_n \leq t_j\}} \mathbb{E} \exp\left\{ i \sum_{k=1}^n \lambda_{\sigma(k)} \cdot X(s_k) \right\} ds_1 \cdots ds_n$$
$$= \sum_{\sigma \in \Sigma_n} \int_{\{0 \leq s_1 \leq \cdots \leq s_n \leq t_j\}} \mathbb{E} \exp\left\{ i \sum_{k=1}^n \left( \sum_{j=k}^n \lambda_{\sigma(j)} \right) \right.$$
$$\left. \times (X(s_k) - X(s_{k-1})) \right\} ds_1 \cdots ds_n$$
$$= \sum_{\sigma \in \Sigma_n} \int_{\{0 \leq s_1 \leq \cdots \leq s_n \leq t_j\}} \prod_{k=1}^n \exp\left\{ -(s_k - s_{k-1}) \psi\left( \sum_{j=k}^n \lambda_{\sigma(j)} \right) \right\} ds_1 \cdots ds_n,$$

where we adopt the convention that $s_0 = 0$. Thus

$$\mathbb{E}[\eta^0([0, t_1] \times \cdots \times [0, t_p])^n]$$

$$(2.1) \quad \begin{aligned} &= \frac{1}{(2\pi)^{dn}} \int_{(\mathbb{R}^d)^n} d\lambda_1 \cdots d\lambda_n \\ &\quad \times \prod_{h=1}^p \sum_{\sigma \in \Sigma_n} \int_{\{0 \le s_1 \le \cdots \le s_n \le t_h\}} \prod_{k=1}^n \exp\bigg\{ -(s_k - s_{k-1}) \\ &\qquad\qquad\qquad\qquad\qquad\qquad \times \psi\bigg( \sum_{j=k}^n \lambda_{\sigma(j)} \bigg) \bigg\} \, ds_1 \cdots ds_n. \end{aligned}$$

To simplify the above representation, we replace $t_1, \ldots, t_p$ by $\tau_1, \ldots, \tau_p$:

$$\begin{aligned} &\mathbb{E}[\eta^0([0, \tau_1] \times \cdots \times [0, \tau_p])^n] \\ &= \frac{1}{(2\pi)^{dn}} \int_{(\mathbb{R}^d)^n} d\lambda_1 \cdots d\lambda_n \\ &\quad \times \bigg[ \sum_{\sigma \in \Sigma_n} \int_0^\infty e^{-t} \, dt \\ &\qquad\qquad \times \int_{\{0 \le s_1 \le \cdots \le s_n \le t\}} \prod_{k=1}^n \exp\bigg\{ -(s_k - s_{k-1}) \\ &\qquad\qquad\qquad\qquad\qquad \times \psi\bigg( \sum_{j=k}^n \lambda_{\sigma(j)} \bigg) \bigg\} \, ds_1 \cdots ds_n \bigg]^p \\ &= \frac{1}{(2\pi)^{dn}} \int_{(\mathbb{R}^d)^n} d\lambda_1 \cdots d\lambda_n \bigg[ \sum_{\sigma \in \Sigma_n} \prod_{k=1}^n \int_0^\infty e^{-t} \exp\bigg\{ -t\psi\bigg( \sum_{j=k}^n \lambda_{\sigma(j)} \bigg) \bigg\} \, dt \bigg]^p \\ &= \frac{1}{(2\pi)^{dn}} \int_{(\mathbb{R}^d)^n} d\lambda_1 \cdots d\lambda_n \bigg[ \sum_{\sigma \in \Sigma_n} \prod_{k=1}^n \bigg[ 1 + \psi\bigg( \sum_{j=k}^n \lambda_{\sigma(j)} \bigg) \bigg]^{-1} \bigg]^p \end{aligned}$$

where the second step follows from the identity ((1.9) in [4]) that

$$\int_0^\infty e^{-t} \, dt \int_{\{0 \le s_1 \le \cdots \le s_n \le t\}} \prod_{k=1}^n \varphi_k(s_k - s_{k-1}) \, ds_1 \cdots ds_n = \prod_{k=1}^n \int_0^\infty e^{-t} \varphi_k(t) \, dt.$$

Write $Q(\lambda) = [1 + \psi(\lambda)]^{-1}$. By the bijection $j \mapsto n - j$ and by the permutation invariance,

$$(2.2) \quad \begin{aligned} &\mathbb{E}[\eta^0([0, \tau_1] \times \cdots \times [0, \tau_p])^n] \\ &= \frac{1}{(2\pi)^{dn}} \int_{(\mathbb{R}^d)^n} d\lambda_1 \cdots d\lambda_n \bigg[ \sum_{\sigma \in \Sigma_n} \prod_{k=1}^n Q\bigg( \sum_{j=1}^k \lambda_{\sigma(j)} \bigg) \bigg]^p. \end{aligned}$$

We state the following theorem which will be proved in Sections 3–5.

THEOREM 2.1. *Under* (1.2),

$$(2.3) \qquad \lim_{n \to \infty} \frac{1}{n} \log \frac{1}{(n!)^p} \mathbb{E}[\eta^0([0, \tau_1] \times \cdots \times [0, \tau_p])^n] = \log \frac{\rho}{(2\pi)^d}$$

*where* $\rho > 0$ *is given in* (1.4).

As it turns out, the hard part of Theorem 1.2 is on the upper bound. On the other hand, if the right constant were not part of our concern, we could establish the upper bound in a much easier way. Indeed, by Jensen's inequality,

$$\int_{(\mathbb{R}^d)^n} d\lambda_1 \cdots d\lambda_n \left[ \sum_{\sigma \in \Sigma_n} \prod_{k=1}^n Q\left( \sum_{j=1}^k \lambda_{\sigma(j)} \right) \right]^p$$

$$\leq (n!)^{p-1} \sum_{\sigma \in \Sigma_n} \int_{(\mathbb{R}^d)^n} d\lambda_1 \cdots d\lambda_n \prod_{k=1}^n Q^p\left( \sum_{j=1}^k \lambda_{\sigma(j)} \right)$$

$$= (n!)^p \int_{(\mathbb{R}^d)^n} d\lambda_1 \cdots d\lambda_n \prod_{k=1}^n Q^p(\lambda_k)$$

$$= (n!)^p \left( \int_{\mathbb{R}^d} Q^p(\lambda) \, d\lambda \right)^n$$

where the second step follows from variable substitution. By (2.2), we obtain the following upper bound:

$$\limsup_{n \to \infty} \frac{1}{n} \log \frac{1}{(n!)^p} \mathbb{E}[\eta^0([0, \tau_1] \times \cdots \times [0, \tau_p])^n] \leq \log\left( \frac{1}{(2\pi)^d} \int_{\mathbb{R}^d} Q^p(\lambda) \, d\lambda \right).$$

Unfortunately, by examining the argument we used to derive (1.5), it is not hard to see that as $p \geq 2$, we would miss the right constant by doing that.

PROOF OF THEOREM 1.1. We now prove Theorem 1.1 based on Theorem 2.1. Let $t_1, \ldots, t_p \geq 0$. In view of (2.1), by Hölder's inequality,

$$\mathbb{E}[\eta^0([0, t_1] \times \cdots \times [0, t_p])^n]$$

$$\leq \frac{1}{(2\pi)^{dn}}$$

$$\times \prod_{h=1}^p \Bigg\{ \int_{(\mathbb{R}^d)^n} d\lambda_1 \cdots d\lambda_n$$

$$\times \Bigg[ \sum_{\sigma \in \Sigma_n} \int_{\{0 \leq s_1 \leq \cdots \leq s_n \leq t_h\}} \prod_{k=1}^n \exp\Bigg\{ -(s_k - s_{k-1})$$

$$\times \psi\left(\sum_{j=k}^{n}\lambda_{\sigma(j)}\right)\bigg\}\,ds_1\cdots ds_n\bigg]^p\bigg\}^{1/p}$$

$$=\prod_{j=1}^{p}\{\mathbb{E}[\eta^0([0,t_j]^p)^n]\}^{1/p}=(t_1\cdots t_p)^{(\alpha p-d)/(\alpha p)n}\mathbb{E}[\eta^0([0,1]^p)^n]$$

where the last step follows from (1.6). Thus,

$$\mathbb{E}[\eta([0,\tau_1]\times\cdots\times[0,\tau_p])^n]$$

$$=\int_0^\infty\cdots\int_0^\infty e^{-(t_1+\cdots+t_p)}\mathbb{E}[\eta([0,t_1]\times\cdots\times[0,t_p])^n]\,dt_1\cdots dt_p$$

$$\leq\mathbb{E}[\eta([0,1]^p)^n]\int_0^\infty\cdots\int_0^\infty (t_1\cdots t_p)^{((\alpha p-d)/(\alpha p))n}e^{-(t_1+\cdots+t_p)}\,dt_1\cdots dt_p$$

$$=\mathbb{E}[\eta([0,1]^p)^n]\left[\Gamma\left(\frac{\alpha p-d}{\alpha p}n+1\right)\right]^p.$$

By Theorem 2.1 and the Stirling formula,

(2.4)
$$\liminf_{n\to\infty}\frac{1}{n}\log(n!)^{-d/\alpha}\mathbb{E}[\eta([0,1]^p)^n]$$

$$\geq\log\left(\frac{\alpha p}{\alpha p-d}\right)^{(\alpha p-d)/\alpha}+\log\frac{\rho}{(2\pi)^d}.$$

On the other hand, notice that $\bar\tau\equiv\min\{\tau_1,\ldots,\tau_p\}$ has the exponential distribution with the parameter $p$. Hence,

$$\mathbb{E}[\eta([0,\tau_1]\times\cdots\times[0,\tau_p])]^n\geq\mathbb{E}[\eta([0,\bar\tau]^p)^n]=\mathbb{E}\bar\tau^{((\alpha p-d)/\alpha)n}\mathbb{E}[\eta([0,1]^p)^n]$$

$$=p^{-((\alpha p-d)/\alpha)n}\Gamma\left(1+\frac{\alpha p-d}{\alpha}n\right)\mathbb{E}[\eta([0,1]^p)^n]$$

where the second step follows from (1.6). By the Stirling formula and Theorem 1.2 we have

(2.5)
$$\limsup_{n\to\infty}\frac{1}{n}\log(n!)^{-d/\alpha}\mathbb{E}[\eta([0,1]^p)^n]$$

$$\leq\log\left(\frac{\alpha p}{\alpha p-d}\right)^{(\alpha p-d)/\alpha}+\log\frac{\rho}{(2\pi)^d}.$$

Combining (2.4) and (2.5) gives

(2.6) $$\lim_{n\to\infty}\frac{1}{n}\log(n!)^{-d/\alpha}\mathbb{E}[\eta([0,1]^p)^n]=\log\left(\frac{\alpha p}{\alpha p-d}\right)^{(\alpha p-d)/\alpha}+\log\frac{\rho}{(2\pi)^d}.$$

Finally, Theorem 1.1 follows from Lemma 1.4. $\square$

**3. Lower bound for Theorem 2.1.** In this section we prove

$$(3.1) \qquad \liminf_{n\to\infty} \frac{1}{n} \log \frac{1}{(n!)^p} \mathbb{E}[\eta^0([0,\tau_1] \times \cdots \times [0,\tau_p])]^n \geq \log \frac{\rho}{(2\pi)^d}.$$

Our starting point is (2.2). Let $q > 1$ be the conjugate number of $p$ defined by $p^{-1} + q^{-1} = 1$ and let $f$ be a positive continuous function on $\mathbb{R}^d$ with $f(-\lambda) = f(\lambda)$ and $\|f\|_q = 1$. We have

$$\left( \int_{(\mathbb{R}^d)^n} d\lambda_1 \cdots d\lambda_n \left[ \sum_{\sigma \in \Sigma_n} \prod_{k=1}^n Q\left( \sum_{j=1}^k \lambda_{\sigma(j)} \right) \right]^p \right)^{1/p}$$

$$\geq \int_{(\mathbb{R}^d)^n} d\lambda_1 \cdots d\lambda_n \left( \prod_{k=1}^n f(\lambda_k) \right) \sum_{\sigma \in \Sigma_n} \prod_{k=1}^n Q\left( \sum_{j=1}^k \lambda_{\sigma(j)} \right)$$

$$= n! \int_{(\mathbb{R}^d)^n} d\lambda_1 \cdots d\lambda_n \left( \prod_{k=1}^n f(\lambda_k) \right) \prod_{k=1}^n Q\left( \sum_{j=1}^k \lambda_j \right)$$

$$= n! \int_{(\mathbb{R}^d)^n} d\lambda_1 \cdots d\lambda_n \prod_{k=1}^n f(\lambda_k - \lambda_{k-1}) Q(\lambda_k)$$

where we follow the convention that $\lambda_0 = 0$.

Define the linear operator $T$ on $\mathcal{L}^2(\mathbb{R}^d)$ as

$$Tg(\lambda) = \sqrt{Q(\lambda)} \int_{\mathbb{R}^d} f(\gamma - \lambda) \sqrt{Q(\gamma)} g(\gamma) \, d\gamma, \qquad g \in \mathcal{L}^2(\mathbb{R}^d).$$

To show that $T$ is well defined and continuous on $\mathcal{L}^2(\mathbb{R}^d)$, we need only to prove that there is a constant $C > 0$ such that

$$(3.2) \qquad \langle h, Tg \rangle \leq C\|g\|_2 \|h\|_2, \qquad g, h \in \mathcal{L}^2(\mathbb{R}^d).$$

Indeed,

$$\langle h, Tg \rangle = \int \int_{\mathbb{R}^d \times \mathbb{R}^d} f(\gamma - \lambda) \sqrt{Q(\lambda)} h(\lambda) \sqrt{Q(\gamma)} g(\gamma) \, d\lambda \, d\gamma$$

$$= \int_{\mathbb{R}^d} f(\gamma) \, d\gamma \int_{\mathbb{R}^d} \sqrt{Q(\lambda)} h(\lambda) \sqrt{Q(\lambda + \gamma)} g(\lambda + \gamma) \, d\lambda$$

$$\leq \left\{ \int_{\mathbb{R}^d} \left[ \int_{\mathbb{R}^d} \sqrt{Q(\lambda)} h(\lambda) \sqrt{Q(\lambda + \gamma)} g(\lambda + \gamma) \, d\lambda \right]^p d\gamma \right\}^{1/p}.$$

Hence, an argument similar to the proof of (1.5) gives that $\langle h, Tg \rangle \leq \|Q\|_p \|g\|_2 \times \|h\|_2$.

In addition, one can see that $\langle h, Tg \rangle = \langle g, Th \rangle$ for any $g, h \in \mathcal{L}^2(\mathbb{R}^d)$. It means that $T$ is self-adjoint. We now let $g$ be a bounded and locally supported function on $\mathbb{R}^d$ with $\|g\|_2 = 1$. Then there is $\delta > 0$ such that $f \geq \delta$ and $Q \geq \delta$ on the support of $g$. In addition, notice that $Q \leq 1$. Thus,

$$\int_{(\mathbb{R}^d)^n} d\lambda_1 \cdots d\lambda_n \prod_{k=1}^{n} f(\lambda_k - \lambda_{k-1}) Q(\lambda_k)$$

$$\geq \delta^2 \|g\|_\infty^{-2} \int_{(\mathbb{R}^d)^n} d\lambda_1 \cdots d\lambda_n \, g(\lambda_1)$$

$$\times \left( \prod_{k=2}^{n} \sqrt{Q(\lambda_{k-1})} f(\lambda_k - \lambda_{k-1}) \sqrt{Q(\lambda_k)} \right) g(\lambda_n)$$

$$= \delta^2 \|g\|_\infty^{-2} \langle g, T^{n-1} g \rangle.$$

Consider the spectral representation of the self-adjoint operator $T$:

$$\langle g, Tg \rangle = \int_{-\infty}^{\infty} \theta \mu_g(d\theta)$$

where $\mu_g(d\theta)$ is a probability measure on $\mathbb{R}$. By the mapping theorem,

$$\langle g, T^{n-1} g \rangle = \int_{-\infty}^{\infty} \theta^{n-1} \mu_g(d\theta) \geq \left( \int_{-\infty}^{\infty} \theta \mu_g(d\theta) \right)^{n-1} = \langle g, Tg \rangle^{n-1}$$

where the second step follows from Jensen's inequality.

Hence,

$$\liminf_{n \to \infty} \frac{1}{n} \log \frac{1}{n!} \left( \int_{(\mathbb{R}^d)^n} d\lambda_1 \cdots d\lambda_n \left[ \sum_{\sigma \in \Sigma_n} \prod_{k=1}^{n} Q\left( \sum_{j=1}^{k} \lambda_{\sigma(j)} \right) \right]^p \right)^{1/p}$$

$$\geq \log \langle g, Tg \rangle = \log \int \int_{\mathbb{R}^d \times \mathbb{R}^d} f(\gamma - \lambda) \sqrt{Q(\lambda)} \sqrt{Q(\gamma)} g(\lambda) g(\gamma) \, d\lambda \, d\gamma$$

$$= \log \int_{\mathbb{R}^d} f(\lambda) \left[ \int_{\mathbb{R}^d} \sqrt{Q(\lambda + \gamma)} \sqrt{Q(\gamma)} g(\lambda + \gamma) g(\gamma) \, d\gamma \right] d\lambda.$$

Notice that the set of all bounded, locally supported $g$ is dense in $\mathcal{L}^2(\mathbb{R}^d)$. Taking the supremum over $g$ on the right-hand sides gives

$$\liminf_{n \to \infty} \frac{1}{n} \log \frac{1}{n!} \left( \int_{(\mathbb{R}^d)^n} d\lambda_1 \cdots d\lambda_n \left[ \sum_{\sigma \in \Sigma_n} \prod_{k=1}^{n} Q\left( \sum_{j=1}^{k} \lambda_{\sigma(j)} \right) \right]^p \right)^{1/p}$$

$$(3.3)$$

$$\geq \log \sup_{|g|_2 = 1} \int_{\mathbb{R}^d} f(\lambda) \left[ \int_{\mathbb{R}^d} \sqrt{Q(\lambda + \gamma)} \sqrt{Q(\gamma)} g(\lambda + \gamma) g(\gamma) \, d\gamma \right] d\lambda$$

Since for any $g$, the function

$$H(\lambda) = \int_{\mathbb{R}^d} \sqrt{Q(\lambda + \gamma)} \sqrt{Q(\gamma)} g(\lambda + \gamma) g(\gamma) \, d\gamma$$

is even: $H(-\lambda) = H(\lambda)$. Hence, taking the supremum over all positive, continuous and even functions $f$ with $\|f\|_q = 1$ on the right-hand side of (3.3) gives

$$\liminf_{n\to\infty} \frac{1}{n} \log \frac{1}{n!} \left( \int_{(\mathbb{R}^d)^n} d\lambda_1 \cdots d\lambda_n \left[ \sum_{\sigma \in \Sigma_n} \prod_{k=1}^n Q\left( \sum_{j=1}^k \lambda_{\sigma(j)} \right) \right]^p \right)^{1/p}$$

$$\geq \frac{1}{p} \log \sup_{|g|_2=1} \int_{\mathbb{R}^d} \left[ \int_{\mathbb{R}^d} \sqrt{Q(\lambda+\gamma)} \sqrt{Q(\gamma)} g(\lambda+\gamma) g(\gamma) \, d\gamma \right]^p d\lambda$$

$$= \frac{1}{p} \log \rho.$$

From the relation (2.2), we have proved (3.1). □

## 4. A discrete version of Theorem 2.1.
The approach for the upper bound of Theorem 2.1 relies heavily on combinatorics and is therefore best suitable for the discrete structure. In this section we prove the following discrete version of Theorem 2.1 with an additional localization assumption.

THEOREM 4.1. *Let $\pi(x)$ and $Q(x)$ be two nonnegative functions on $\mathbb{Z}^d$ such that $\pi$ is locally supported, $\pi(-x) = \pi(x)$ for all $x \in \mathbb{Z}^d$, and that*

$$(4.1) \qquad\qquad \lim_{|x|\to\infty} Q(x) = 0.$$

*Then*

$$\lim_{n\to\infty} \frac{1}{n} \log \sum_{x_1,\ldots,x_n \in \mathbb{Z}^d} \left( \prod_{k=1}^n \pi(x_k) \right) \left[ \frac{1}{n!} \sum_{\sigma \in \Sigma_n} \prod_{k=1}^n Q\left( \sum_{j=1}^k x_{\sigma(j)} \right) \right]^p$$

$$(4.2)$$

$$= \log \tilde{\rho},$$

*where*

$$\tilde{\rho} = \sup_{|f|_2=1} \sum_{x \in \mathbb{Z}^d} \pi(x) \left[ \sum_{y \in \mathbb{Z}^d} \sqrt{Q(x+y)} \sqrt{Q(y)} f(x+y) f(y) \right]^p$$

*and*

$$|f|_2 = \left( \sum_{x \in \mathbb{Z}^d} f^2(x) \right)^{1/2}.$$

PROOF. The lower bound follows from an obvious modification of the argument in the previous section. We now prove the upper bound. By assumption, there is a finite set $A \subset \mathbb{Z}^d$ such that $\pi(x) > 0$ as $x \in A$ and $\pi(x) = 0$ as $x \notin A$. Without loss of generality, we may assume that the group generated by $A$ is $\mathbb{Z}^d$. Indeed, if $A$ does not generate $\mathbb{Z}^d$, one can add finitely many lattice points into $A$ to form an augmented $\bar{A}$ which generates $\mathbb{Z}^d$. Let $\varepsilon > 0$ be a small number. Assume that we have proved the upper bound under this extra condition. We apply it to the system where $\pi(\cdot)$ is replaced by $\bar{\pi}(\cdot)$ defined as: $\bar{\pi}(x) = \pi(x)$ on $A \cup (\mathbb{Z}^d \setminus \bar{A})$ and $\bar{\pi}(x) = \varepsilon$ on $\bar{A} \setminus A$:

$$\limsup_{n \to \infty} \frac{1}{n} \log \sum_{x_1,\ldots,x_n \in \mathbb{Z}^d} \left( \prod_{k=1}^n \pi(x_k) \right) \left[ \frac{1}{n!} \sum_{\sigma \in \Sigma_n} \prod_{k=1}^n Q\left( \sum_{j=1}^k x_{\sigma(j)} \right) \right]^p$$

$$\leq \limsup_{n \to \infty} \frac{1}{n} \log \sum_{x_1,\ldots,x_n \in \mathbb{Z}^d} \left( \prod_{k=1}^n \bar{\pi}(x_k) \right) \left[ \frac{1}{n!} \sum_{\sigma \in \Sigma_n} \prod_{k=1}^n Q\left( \sum_{j=1}^k x_{\sigma(j)} \right) \right]^p$$

$$\leq \log \sup_{\|f\|_2 = 1} \sum_{x \in \mathbb{Z}^d} \bar{\pi}(x) \left[ \sum_{y \in \mathbb{Z}^d} \sqrt{Q(x+y)} \sqrt{Q(y)} f(x+y) f(y) \right]^p.$$

Letting $\varepsilon \to 0^+$ on the right-hand side gives the desired upper bound.

We may also assume that $\pi$ is a probability measure on $A$, for otherwise we use $\pi(\cdot)/\pi(A)$ instead of $\pi(\cdot)$ in the following proof.

We adopt the notation $\mathbf{y} = (y_1, \ldots, y_n)$ for any $y_1, \ldots, y_n \in \mathbb{Z}^d$ and write

$$L_n^{\mathbf{y}} = \frac{1}{n} \sum_{k=1}^n \delta_{y_k}.$$

Let $n$ and $\mathbf{x} = (x_1, \ldots, x_n) \in A^n$ be fixed for a moment and write $\mu = L_n^{\mathbf{x}}$. Then for each $x \in A$, $n\mu(x)$ is an integer, and

$$\sum_{\sigma \in \Sigma_n} \prod_{k=1}^n Q\left( \sum_{j=1}^k x_{\sigma(j)} \right)$$

$$= \sum_{y_1,\ldots,y_n \in A} \mathbb{1}_{\{L_n^{\mathbf{y}} = \mu\}} \sum_{\sigma \in \Sigma_n} \mathbb{1}_{\{\mathbf{x} \circ \sigma = \mathbf{y}\}} \prod_{k=1}^n Q\left( \sum_{j=1}^k y_j \right)$$

$$= \sum_{y_1,\ldots,y_n \in A} \mathbb{1}_{\{L_n^{\mathbf{y}} = \mu\}} \left\{ \prod_{k=1}^n Q\left( \sum_{j=1}^k y_j \right) \right\} \#\{\sigma \in \Sigma_n;\ \mathbf{x} \circ \sigma = \mathbf{y}\}.$$

Notice that as $L_n^{\mathbf{y}} = \mu$,

$$(4.3) \qquad \#\{\sigma \in \Sigma_n;\ \mathbf{x} \circ \sigma = \mathbf{y}\} = \prod_{x \in A} (n\mu(x))!.$$

Indeed, for each $x \in A$ there are, respectively, exactly $n\mu(x)$ of $x_1, \ldots, x_n$ and exactly $n\mu(x)$ of $y_1, \ldots, y_n$ which are equal to $x$. Therefore, there are $(n\mu(x))!$ ways to match each $x$-valued component of $\mathbf{y}$ to each $x$-valued component of $\mathbf{x}$. Thus, (4.3) follows from the multiplication principle.

Hence

$$
\begin{aligned}
(4.4) \qquad & \sum_{\sigma \in \Sigma_n} \prod_{k=1}^{n} Q\left(\sum_{j=1}^{k} x_{\sigma(j)}\right) \\
& = \left(\prod_{x \in A} (n\mu(x))!\right) \sum_{y_1, \ldots, y_n \in A} \mathbb{1}_{\{L_n^{\mathbf{y}} = \mu\}} \left\{\prod_{k=1}^{n} Q\left(\sum_{j=1}^{k} y_j\right)\right\}.
\end{aligned}
$$

By the Stirling formula, $n! \sim \sqrt{2\pi n} n^n e^{-n}$ $(n \to \infty)$ and there is $C > 0$ such that

$$
(n\mu(x))! \le C\sqrt{n\mu(x)}(n\mu(x))^{n\mu(x)} e^{-n\mu(x)}
$$

for all $x \in A$ and all $n \ge 1$. Consequently,

$$
\begin{aligned}
\prod_{x \in A} (n\mu(x))! & \le C n^{\#(A)/2} e^{-n} \prod_{x \in A} n^{n\mu(x)} \exp\{n\mu(x) \log \mu(x)\} \\
& = C n^{\#(A)/2} e^{-n} n^n \exp\left\{n \sum_{x \in A} \mu(x) \log \mu(x)\right\} \\
& \le C n^{\#(A)/2} n! \exp\left\{n \sum_{x \in A} \mu(x) \log \mu(x)\right\}.
\end{aligned}
$$

Therefore,

$$
\begin{aligned}
(4.5) \qquad \frac{1}{n!} \sum_{\sigma \in \Sigma_n} \prod_{k=1}^{n} Q\left(\sum_{j=1}^{k} x_{\sigma(j)}\right) & \le C n^{\#(A)/2} \exp\left\{n \sum_{x \in A} \mu(x) \log \mu(x)\right\} \\
& \times \sum_{y_1, \ldots, y_n \in A} \mathbb{1}_{\{L_n^{\mathbf{y}} = \mu\}} \left\{\prod_{k=1}^{n} Q\left(\sum_{j=1}^{k} y_j\right)\right\}.
\end{aligned}
$$

Here and elsewhere below, we follow the convention $0^0 = 1$ or, $0 \log 0 = 0$.

On the other hand, let $q > 1$ be the conjugate number of $p$ defined by $p^{-1} + q^{-1} = 1$. For any probability measure $\nu$ on $A$, write

$$
\phi_\nu(x) = \begin{cases} (\nu(x))^{1/q}(\pi(x))^{1/p}, & x \in A, \\ 0, & x \in \mathbb{Z}^d \setminus A. \end{cases}
$$

Notice that as $L_n^{\mathbf{y}} = \mu$, there are exactly $n\mu(x)$ of $\phi_\mu(y_1), \ldots, \phi_\mu(y_n)$ equal to $\phi_\mu(x)$ for each $x \in A$. Hence,

$$
\sum_{y_1, \ldots, y_n \in A} \phi_\mu(y_1) \cdots \phi_\mu(y_n) \left\{\prod_{k=1}^{n} Q\left(\sum_{j=1}^{k} y_j\right)\right\}
$$

$$\geq \sum_{y_1,\ldots,y_n \in A} \mathbb{1}_{\{L_n^{\mathbf{y}}=\mu\}} \phi_\mu(y_1) \cdots \phi_\mu(y_n) \left\{ \prod_{k=1}^{n} Q\left(\sum_{j=1}^{k} y_j\right) \right\}$$

$$= \sum_{y_1,\ldots,y_n \in A} \mathbb{1}_{\{L_n^{\mathbf{y}}=\mu\}} \left( \prod_{x \in A} \phi_\mu(x)^{n\mu(x)} \right) \prod_{k=1}^{n} Q\left(\sum_{j=1}^{k} y_j\right)$$

$$= \exp\left\{ n \left( \frac{1}{q} \sum_{x \in A} \mu(x) \log \mu(x) + \frac{1}{p} \sum_{x \in A} \mu(x) \log \pi(x) \right) \right\}$$

$$\times \sum_{y_1,\ldots,y_n \in A} \mathbb{1}_{\{L_n^{\mathbf{y}}=\mu\}} \left\{ \prod_{k=1}^{n} Q\left(\sum_{j=1}^{k} y_j\right) \right\}.$$

Combining this with (4.5),

$$\frac{1}{n!} \sum_{\sigma \in \Sigma_n} \prod_{k=1}^{n} Q\left(\sum_{j=1}^{k} x_{\sigma(j)}\right)$$

$$\leq C n^{\#(A)/2} \exp\left\{ n\frac{1}{p} \sum_{x \in A} \mu(x) \log \frac{\mu(x)}{\pi(x)} \right\}$$

$$\times \sum_{y_1,\ldots,y_n \in \mathbb{Z}^d} \phi_\mu(y_1) \cdots \phi_\mu(y_n) \left\{ \prod_{k=1}^{n} Q\left(\sum_{j=1}^{k} y_j\right) \right\}.$$

By variable substitution,

$$\sum_{y_1,\ldots,y_n \in \mathbb{Z}^d} \phi_\mu(y_1) \cdots \phi_\mu(y_n) \left\{ \prod_{k=1}^{n} Q\left(\sum_{j=1}^{k} y_j\right) \right\}$$

$$= \sum_{z_1,\ldots,z_n \in \mathbb{Z}^d} \prod_{k=1}^{n} \phi_\mu(z_k - z_{k-1}) Q(z_k).$$

Summarizing what we have proved,

$$\sum_{x_1,\ldots,x_n \in \mathbb{Z}^d} \left( \prod_{k=1}^{n} \pi(x_k) \right) \left[ \frac{1}{n!} \sum_{\sigma \in \Sigma_n} \prod_{k=1}^{n} Q\left(\sum_{j=1}^{k} x_{\sigma(j)}\right) \right]^p$$

$$(4.6) \qquad \leq C n^{p\#(A)/2} \left[ \sup_{\nu \in \mathcal{M}(A)} \sum_{z_1,\ldots,z_n \in \mathbb{Z}^d} \prod_{k=1}^{n} \phi_\nu(z_k - z_{k-1}) Q(z_k) \right]^p$$

$$\times \sum_{x_1,\ldots,x_n \in A} \pi(x_1) \cdots \pi(x_n) \exp\left\{ n \sum_{x \in A} L_n^{\mathbf{x}}(x) \log \frac{L_n^{\mathbf{x}}(x)}{\pi(x)} \right\},$$

where $\mathcal{M}(A)$ is the space of all probability measures on $A$ equipped with topology of weak convergence. (In our setting, of course, the weak convergence is equivalent to the pointwise convergence.) Recall that by Sanov's theorem (Theorem 2.1.10, page 16 in [12]), the empirical measure $L_n^\mathbf{x}$ satisfies the large deviation principle governed by the rate function

$$H(\nu|\pi) = \sum_{x \in A} \nu(x) \log \frac{\nu(x)}{\pi(x)}, \qquad \nu \in \mathcal{M}(A).$$

By the fact that $A$ is finite and that $\pi(x) > 0$ on $A$, $H(\nu|\pi)$ is continuous on $\mathcal{M}(A)$. By Varadhan's integral lemma (Theorem 4.3.1, page 137 in [12]),

$$\lim_{n \to \infty} \frac{1}{n} \log \sum_{x_1, \ldots, x_n \in A} \pi(x_1) \cdots \pi(x_n) \exp\left\{ n \sum_{x \in A} L_n^\mathbf{x}(x) \log \frac{L_n^\mathbf{x}(x)}{\pi(x)} \right\}$$
$$= \sup_{\nu \in \mathcal{M}(A)} \{ H(\nu|\pi) - H(\nu|\pi) \} = 0.$$

In view of (4.6), the conclusion follows from the following Lemma 4.2. □

LEMMA 4.2. *Under the assumptions given above,*

$$\limsup_{n \to \infty} \frac{1}{n} \log \sup_{\nu \in \mathcal{M}(A)} \sum_{x_1, \ldots, x_n \in \mathbb{Z}^d} \prod_{k=1}^n \phi_\nu(x_k - x_{k-1}) Q(x_k) \le \frac{1}{p} \log \tilde\rho.$$

PROOF. Notice that $\mathcal{M}(A)$ is a compact space and that for any $\mu_0 \in \mathcal{M}(A)$ and $\varepsilon > 0$, there is a open neighborhood $\mathcal{U}$ of $\mu_0$ such that $\mu(x) \le \mu_0(x) + \varepsilon$ for all $\mu \in \mathcal{U}$. Fix $\mu_0$ and write

$$\varphi_\varepsilon(x) = \begin{cases} (\mu_0(x) + \varepsilon)^{1/q} (\pi(x))^{1/p}, & x \in A, \\ 0, & x \in \mathbb{Z}^d \setminus A, \end{cases}$$

and

$$\Lambda_\varepsilon(\mu_0) = \limsup_{n \to \infty} \frac{1}{n} \log \sum_{x_1, \ldots, x_n \in \mathbb{Z}^d} \prod_{k=1}^n \varphi_\varepsilon(x_k - x_{k-1}) Q(x_k).$$

We need only to show that

$$(4.7) \qquad \limsup_{\varepsilon \to 0^+} \Lambda_\varepsilon(\mu_0) \le \frac{1}{p} \log \tilde\rho \qquad \text{uniformly over } \mu_0 \in \mathcal{M}(A).$$

For any $x = (x_1, \ldots, x_d) \in \mathbb{Z}^d$ write $|x|_\infty = \max_{1 \le i \le d} |x_i|$. For any integer $a < b$ we use $(a, b)^d$ and $[a, b]^d$ below for the $d$-dimensional boxes of lattice points. Let $N_0 = \max\{|x|_\infty;\ x \in A\}$. Let $\delta > 0$ be fixed and take integer $N > 2N_0$ sufficiently large so that $Q(x) \leq \delta$ for all $x \in \mathbb{Z}^d$ with $|x|_\infty \geq N/2$. We have $A \subset (-N, N]^d$.

$$
\begin{aligned}
\sum_{x_1,\ldots,x_n \in \mathbb{Z}^d} &\prod_{k=1}^{n} \varphi_\varepsilon(x_k - x_{k-1}) Q(x_k) \\
&= \sum_{y_1,\ldots,y_n \in \mathbb{Z}^d} \sum_{z_1,\ldots,z_n \in (-N,N]^d} \prod_{k=1}^{n} \varphi_\varepsilon(2(y_k - y_{k-1})N \\
&\hspace{4cm} + (z_k - z_{k-1}))Q(2y_k N + z_k) \\
&\leq \sum_{z_1,\ldots,z_n \in (-N,N]^d} \prod_{k=1}^{n} \widetilde{\varphi}_\varepsilon(z_k - z_{k-1}) Q^*(z_k),
\end{aligned}
\tag{4.8}
$$

where

$$
\widetilde{\varphi}_\varepsilon(x) = \sum_{y \in \mathbb{Z}^d} \varphi_\varepsilon(2yN + x), \qquad Q^*(x) = \sup_{y \in \mathbb{Z}^d} Q(2yN + x).
$$

We have

$$
\widetilde{\varphi}_\varepsilon(x) = \varphi_\varepsilon(x), \qquad x \in [-(2N - N_0), (2N - N_0)]^d,
\tag{4.9}
$$

$$
\sum_{y \in (-N,N]^d} \widetilde{\varphi}_\varepsilon(y - x) = \sum_{y \in (-N,N]^d} \widetilde{\varphi}_\varepsilon(y), \qquad x \in (-N, N]^d,
\tag{4.10}
$$

$$
Q^*(x) \leq \delta \vee Q(x), \qquad x \in \mathbb{Z}^d.
\tag{4.11}
$$

By (4.10), the kernel

$$
P(x, y) = \widetilde{\varphi}_\varepsilon(y - x) \Big/ \sum_{z \in (-N,N]^d} \widetilde{\varphi}_\varepsilon(z), \qquad x, y \in (-N, N]^d,
$$

is a transition probability on $(-N, N]^d$. Let $\{Y_k\}_{k \geq 1}$ be a Markov chain with the transition $P(x, y)$. By (4.9), by the definition of $\varphi_\varepsilon$ and by the assumption that the group generated by $A$ is $\mathbb{Z}^d$, $\{Y_k\}_{k \geq 1}$ is irreducible. By the large deviation principle for the empirical measures of finite-state Markov chains (Theorems 3.1.2 and 3.1.6 in [12]), the empirical measure

$$
L_n^{\mathbf{Y}} = \frac{1}{n} \sum_{k=1}^{n} \delta_{Y_k}
$$

satisfies the large deviation principle on $\mathcal{M}\{(-N, N]^d\}$ governed by the rate function

$$
\begin{aligned}
I(\mu) = -\inf_{u > 0} \sum_{x \in (-N,N]^d} \mu(x) \log\left( u(x)^{-1} \sum_{y \in (-N,N]^d} P(x, y) u(y) \right), \\
\mu \in \mathcal{M}\{(-N, N]^d\}.
\end{aligned}
$$

On the other hand,

$$\sum_{z_1,\ldots,z_n \in (-N,N]^d} \prod_{k=1}^{n} \widetilde{\varphi}_\varepsilon(z_k - z_{k-1})Q^*(z_k)$$

$$= \left(\sum_{x \in (-N,N]^d} \widetilde{\varphi}_\varepsilon(x)\right)^n \mathbb{E}_0 \exp\{n\langle \log Q^*, L_n^{\mathbf{Y}}\rangle\}.$$

By Varadhan's integral lemma (Theorem 4.3.1, page 137 in [12]),

$$\lim_{n\to\infty} \frac{1}{n} \log \sum_{z_1,\ldots,z_n \in (-N,N]^d} \prod_{k=1}^{n} \widetilde{\varphi}_\varepsilon(z_k - z_{k-1})Q^*(z_k)$$

$$= \log\left(\sum_{x \in (-N,N]^d} \widetilde{\varphi}_\varepsilon(x)\right)$$

$$+ \sup_{\mu \in \mathcal{M}\{(-N,N]^d\}} \left\{\sum_{x \in (-N,N]^d} \mu(x) \log Q^*(x) - I(\mu)\right\}$$

$$(4.12) \quad = \sup_{\mu \in \mathcal{M}\{(-N,N]^d\}} \inf_{u>0} \sum_{x \in (-N,N]^d} \mu(x) \log\Bigg(Q^*(x)u(x)^{-1}$$

$$\times \sum_{y \in (-N,N]^d} \widetilde{\varphi}_\varepsilon(y-x)u(y)\Bigg)$$

$$\leq \sup_{\mu \in \mathcal{M}\{(-N,N]^d\}} \inf_{u>0} \log\left(\sum_{x,y \in (-N,N]^d} \mu(x)Q^*(x)u(x)^{-1}\widetilde{\varphi}_\varepsilon(y-x)u(y)\right)$$

where the last step follows from Jensen's inequality.

Let $u(x) = \sqrt{Q^*(x)\mu(x)}$ and

$$f(x) = \begin{cases} \sqrt{\mu(x)}, & x \in (-N,N]^d, \\ 0, & x \in \mathbb{Z}^d \setminus (-N,N]^d. \end{cases}$$

We have $|f|_2 = 1$ and

$$\sum_{x,y \in (-N,N]^d} \mu(x)Q^*(x)u(x)^{-1}\widetilde{\varphi}_\varepsilon(y-x)u(y)$$

$$= \sum_{x,y \in (-N,N]^d} \widetilde{\varphi}_\varepsilon(y-x)\sqrt{Q^*(x)Q^*(y)}f(x)f(y).$$

By (4.9), for any $x, y \in (-N, N]^d$, $x - y \notin [-(2N - N_0), (2N - N_0)]^d$ implies that $|x|_\infty \geq N - N_0$ and $|y|_\infty \geq N - N_0$. In view of (4.11),

$$\sum_{x,y \in (-N,N]^d} \widetilde{\varphi}_\varepsilon(y - x) \sqrt{Q^*(x) Q^*(y)} f(x) f(y)$$

$$\leq \sum_{x,y \in \mathbb{Z}^d} \varphi_\varepsilon(y - x) \sqrt{Q_\delta(x) Q_\delta(y)} f(x) f(y)$$

$$+ \sum_{x,y \in B} \widetilde{\varphi}_\varepsilon(y - x) \sqrt{Q^*(x) Q^*(y)} f(x) f(y),$$

where $Q_\delta(x) = \delta \vee Q(x)$ and $B = \{x \in (-N, N]^d; |x|_\infty \geq N - N_0\}$. By the fact [partially from (4.11)] that $Q^*(x) \leq \delta$ for $x \in B$,

$$\sum_{x,y \in B} \widetilde{\varphi}_\varepsilon(y - x) \sqrt{Q^*(x) Q^*(y)} f(x) f(y) \leq \delta \sum_{x,y \in (-N,N]^d} \widetilde{\varphi}_\varepsilon(y - x) f(x) f(y).$$

To control the right-hand side, we consider Fourier transformation. For any function $g$ supported on $(-N, N]^d$, we introduce the complex function $\mathcal{F}(g)$ on $\mathbb{Z}^d$ by

$$\mathcal{F}(g)(y) = \sum_{x \in (-N,N]^d} g(x) \exp\left\{ i \frac{\pi}{N} (x \cdot y) \right\}, \qquad y \in \mathbb{Z}^d.$$

By orthogonality, for any $g$ and $h$ supported on $(-N, N]^d$,

$$(4.13) \qquad \sum_{y \in (-N,N]^d} \mathcal{F}(g)(y) \overline{\mathcal{F}(h)(y)} = (2N)^d \sum_{x \in (-N,N]^d} g(x) h(x).$$

We now take

$$h(x) = \sum_{z \in (-N,N]^d} \widetilde{\varphi}_\varepsilon(z - x) f(z).$$

Then

$$\mathcal{F}(h)(y) = \sum_{z \in (-N,N]^d} f(z) \sum_{x \in (-N,N]^d} \widetilde{\varphi}_\varepsilon(z - x) \exp\left\{ i \frac{\pi}{N} (x \cdot y) \right\}$$

$$= \sum_{z \in (-N,N]^d} f(z) \exp\left\{ i \frac{\pi}{N} (z \cdot y) \right\}$$

$$\times \sum_{x \in (-N,N]^d} \widetilde{\varphi}_\varepsilon(z - x) \exp\left\{ i \frac{\pi}{N} ((x - z) \cdot y) \right\}$$

$$= \mathcal{F}(f)(y) \mathcal{F}(\widetilde{\varphi}_\varepsilon)(-y),$$

where the last step partially follows from the fact that $\widetilde{\varphi}_\varepsilon$ is periodic:

$$\widetilde{\varphi}_\varepsilon(x + 2Ny) = \widetilde{\varphi}_\varepsilon(x), \qquad x \in (-N, N]^d, \ y \in \mathbb{Z}^d.$$

From (4.13), therefore,

$$\sum_{x,y \in (-N,N]^d} \widetilde{\varphi}_\varepsilon(y-x) f(x) f(y)$$
$$= \sum_{x \in (-N,N]^d} f(x) h(x)$$
$$= (2N)^{-d} \sum_{y \in (-N,N]^d} \mathcal{F}(f)(y) \overline{\mathcal{F}(h)(y)}$$
$$= (2N)^{-d} \sum_{y \in (-N,N]^d} |\mathcal{F}(f)(y)|^2 \mathcal{F}(\widetilde{\varphi}_\varepsilon)(y).$$

By the definition of $\widetilde{\varphi}_\varepsilon$,

$$\mathcal{F}(\widetilde{\varphi}_\varepsilon)(y) = \sum_{z \in \mathbb{Z}^d} \sum_{x \in (-N,N]^d} \varphi_\varepsilon(2Nz + x) \exp\left\{ i \frac{\pi}{N}(x \cdot y) \right\}$$
$$= \sum_{z \in \mathbb{Z}^d} \sum_{x \in (-N,N]^d} \varphi_\varepsilon(2Nz + x) \exp\left\{ i \frac{\pi}{N}((2Nz + x) \cdot y) \right\}$$
$$= \sum_{x \in \mathbb{Z}^d} \varphi_\varepsilon(x) \exp\left\{ i \frac{\pi}{N}(x \cdot y) \right\} = \sum_{x \in A} \varphi_\varepsilon(x) \exp\left\{ i \frac{\pi}{N}(x \cdot y) \right\}.$$

Thus, there is a constant $C > 0$ independent of $N$ (and therefore $\delta$), such that $|\mathcal{F}(\widetilde{\varphi}_\varepsilon)(y)| \le C$ for any $y \in (-N, N]^d$.

Therefore, by (4.13) again,

$$\sum_{x,y \in (-N,N]^d} \widetilde{\varphi}_\varepsilon(y-x) f(x) f(y) \le C(2N)^{-d} \sum_{y \in (-N,N]^d} |\mathcal{F}(f)(y)|^2$$
$$= C \sum_{x \in (-N,N]^d} f^2(x) = C.$$

Summarizing the above discussion, by (4.12) we have

$$\lim_{n \to \infty} \frac{1}{n} \log \sum_{z_1, \ldots, z_n \in (-N,N]^d} \prod_{k=1}^n \widetilde{\varphi}_\varepsilon(z_k - z_{k-1}) Q^*(z_k)$$
$$\le \log\left( C\delta + \sup_{|f|_2 = 1} \sum_{x,y \in \mathbb{Z}^d} \varphi_\varepsilon(y - x) \sqrt{Q_\delta(x) Q_\delta(y)} f(x) f(y) \right).$$

By (4.8),

$$\limsup_{n\to\infty} \frac{1}{n} \log \sum_{x_1,\ldots,x_n\in\mathbb{Z}^d} \prod_{k=1}^n \varphi_\varepsilon(x_k-x_{k-1})Q(x_k)$$

$$\leq \log\left(C\delta + \sup_{|f|_2=1} \sum_{x,y\in\mathbb{Z}^d} \varphi_\varepsilon(y-x)\sqrt{Q_\delta(x)Q_\delta(y)}f(x)f(y)\right).$$

Letting $\delta\to 0^+$ gives

$$\limsup_{n\to\infty} \frac{1}{n} \log \sum_{x_1,\ldots,x_n\in\mathbb{Z}^d} \prod_{k=1}^n \varphi_\varepsilon(x_k-x_{k-1})Q(x_k)$$

$$\leq \log\left(\sup_{|f|_2=1} \sum_{x,y\in\mathbb{Z}^d} \varphi_\varepsilon(y-x)\sqrt{Q(x)Q(y)}f(x)f(y)\right).$$

By the definition of $\varphi_\varepsilon$, for any $f\in\mathcal{L}^2(\mathbb{Z}^d)$ with $|f|_2=1$,

$$\sum_{x,y\in\mathbb{Z}^d} \varphi_\varepsilon(y-x)\sqrt{Q(x)Q(y)}f(x)f(y)$$

$$= \sum_{x\in\mathbb{Z}^d} \varphi_\varepsilon(x) \sum_{y\in\mathbb{Z}^d} \sqrt{Q(x+y)Q(y)}f(x+y)f(y)$$

$$= \sum_{x\in A} (\mu_0(x)+\varepsilon)^{1/q}(\pi(x))^{1/p} \sum_{y\in\mathbb{Z}^d} \sqrt{Q(x+y)Q(y)}f(x+y)f(y)$$

$$\leq \left\{\sum_{x\in A} (\mu_0(x)+\varepsilon)\right\}^{1/q}$$

$$\times \left\{\sum_{x\in A} \pi(x)\left[\sum_{y\in\mathbb{Z}^d} \sqrt{Q(x+y)Q(y)}f(x+y)f(y)\right]^p\right\}^{1/p}$$

$$\leq (1+\varepsilon\#\{A\})^{1/q}\tilde{\rho}^{1/p}.$$

Consequently,

$$\limsup_{n\to\infty} \frac{1}{n} \log \sum_{x_1,\ldots,x_n\in\mathbb{Z}^d} \prod_{k=1}^n \varphi_\varepsilon(x_k-x_{k-1})Q(x_k)$$

$$\leq \frac{1}{p}\log\tilde{\rho} + \frac{1}{q}\log(1+\varepsilon\#\{A\})$$

which clearly implies (4.7). $\square$

**5. Upper bound for Theorem 2.1.** In this section we prove

$$(5.1) \qquad \limsup_{n \to \infty} \frac{1}{n} \log \frac{1}{(n!)^p} \mathbb{E}[\eta^0([0, \tau_1] \times \cdots \times [0, \tau_p])^n] \leq \log \frac{\rho}{(2\pi)^d}.$$

By comparing (2.2) with Theorem 4.1, we need to do two things—localization and discretization. In particular, we point out the difficulty in our second task. If we follow a standard way of discretization, then each of $\lambda_1, \ldots, \lambda_n$ will generate a small error. This may lead to a considerable error generated by

$$\sum_{j=1}^{k} \lambda_{\sigma(j)}$$

as $k$ is large. In view of (2.2), therefore, the standard approach seems not to be very promising.

Our approach relies on Fourier transformation. Define the probability density $h$ on $\mathbb{R}^d$ as

$$(5.2) \qquad h(x) = C^{-1} \prod_{k=1}^{d} \left( \frac{2 \sin x_k}{x_k} \right)^2, \qquad x = (x_1, \ldots, x_d) \in \mathbb{R}^d,$$

where $C > 0$ is the normalizing constant:

$$C = \int_{\mathbb{R}^d} \prod_{k=1}^{d} \left( \frac{2 \sin x_k}{x_k} \right)^2 dx_1 \cdots dx_d.$$

Clearly, $h$ is symmetric. One can verify that the Fourier transform $\widehat{h}$ is

$$\widehat{h}(\lambda) = \int_{\mathbb{R}^d} h(x) e^{i\lambda \cdot x} \, dx = C^{-1} (2\pi)^d (\mathbb{1}_{[-1,1]^d} * \mathbb{1}_{[-1,1]^d})(\lambda).$$

In particular, $\widehat{h}$ is nonnegative and has the compact support set $[-2, 2]^d$.

For each $\varepsilon > 0$, write

$$(5.3) \qquad h_\varepsilon(x) = \varepsilon^{-d} h(\varepsilon^{-1} x), \qquad x \in \mathbb{R}^d.$$

Define

$$\eta_\varepsilon([0, t_1] \times \cdots \times [0, t_p]) = \int_{\mathbb{R}^d} h_\varepsilon(x) \eta^x([0, t_1] \times \cdots \times [0, t_p]) \, dx.$$

By Parseval's identity we have

$$\eta_\varepsilon([0, t_1] \times \cdots \times [0, t_p])$$
$$= \frac{1}{(2\pi)^d} \int_{\mathbb{R}^d} d\lambda \, \widehat{h}_\varepsilon(\lambda)$$
$$\times \int_0^{t_1} \cdots \int_0^{t_p} \exp\{i\lambda \cdot (X_1(s_1) + \cdots + X_p(s_p))\} \, ds_1 \cdots ds_p,$$

where

$$\widehat{h}_\varepsilon(\lambda) = \int_{\mathbb{R}^d} h_\varepsilon(x) e^{i\lambda \cdot x} \, dx = \widehat{h}(\varepsilon\lambda).$$

Hence,

$$\mathbb{E}[(\eta^0 - \eta_\varepsilon)([0,t_1] \times \cdots \times [0,t_p])^n]$$

$$= \frac{1}{(2\pi)^{nd}} \int_{(\mathbb{R}^d)^n} d\lambda_1 \cdots d\lambda_n \left( \prod_{k=1}^n [1 - \widehat{h}(\varepsilon\lambda_k)] \right)$$

$$\times \prod_{j=1}^p \int_{[0,t_j]^n} \mathbb{E} \exp\left\{ i \sum_{k=1}^n \lambda_k \cdot X(s_k) \right\} ds_1 \cdots ds_n.$$

Following the same procedure used for (2.2),

$$\mathbb{E}[(\eta^0 - \eta_\varepsilon)([0,\tau_1] \times \cdots \times [0,\tau_p])^n]$$

$$= \frac{1}{(2\pi)^{dn}} \int_{(\mathbb{R}^d)^n} d\lambda_1 \cdots d\lambda_n \left( \prod_{k=1}^n [1 - \widehat{h}(\varepsilon\lambda_k)] \right)$$

(5.4)

$$\times \left[ \sum_{\sigma \in \Sigma_n} \prod_{k=1}^n Q\left( \sum_{j=1}^k \lambda_{\sigma(j)} \right) \right]^p$$

$$\leq \frac{(n!)^p}{(2\pi)^{dn}} \int_{(\mathbb{R}^d)^n} d\lambda_1 \cdots d\lambda_n \prod_{k=1}^n [1 - \widehat{h}(\varepsilon\lambda_k)] Q^p\left( \sum_{j=1}^k \lambda_j \right)$$

$$= \frac{(n!)^p}{(2\pi)^{dn}} \int_{(\mathbb{R}^d)^n} d\lambda_1 \cdots d\lambda_n \prod_{k=1}^n [1 - \widehat{h}(\varepsilon(\lambda_k - \lambda_{k-1}))] Q^p(\lambda_k),$$

where $Q(\lambda) = [1 + \psi(\lambda)]^{-1}$, where the second step follows from Hölder inequality and from a suitable index rearrangement, and where the third step follows from the variable substitution $\lambda_k \mapsto \lambda_k - \lambda_{k-1}$ (recall our convention $\lambda_0 = 0$).

We now prove that

(5.5)

$$\limsup_{\varepsilon \to 0^+} \limsup_{n \to \infty} \frac{1}{n} \log \int_{(\mathbb{R}^d)^n} d\lambda_1 \cdots d\lambda_n$$

$$\times \prod_{k=1}^n [1 - \widehat{h}(\varepsilon(\lambda_k - \lambda_{k-1}))] Q^p(\lambda_k) = -\infty.$$

First notice that under the assumption $d < \alpha p$,

$$C \equiv \int_{\mathbb{R}^d} Q^p(\lambda) \, d\lambda < \infty.$$

Given $\delta > 0$ there are $u > 0$ and $N > 0$ such that $1 - \widehat{h}(\lambda) < \delta$ as $|\lambda| < u$, and that

$$\int_{\{|\lambda| \geq N\}} Q^p(\lambda) \, d\lambda < \delta.$$

We take $\varepsilon < u(2N)^{-1}$. For each $n$, write

$$B_n = \Big\{ (\lambda_1, \dots, \lambda_n) \in (\mathbb{R}^d)^n; \ \#\{1 \leq k \leq n; \ |\lambda_k| \geq N\} \geq \frac{n}{3} \Big\}.$$

We have

$$\int_{(\mathbb{R}^d)^n} d\lambda_1 \cdots d\lambda_n \prod_{k=1}^n [1 - \widehat{h}(\varepsilon(\lambda_k - \lambda_{k-1}))] Q^p(\lambda_k)$$

$$\leq \int_{B_n} d\lambda_1 \cdots d\lambda_n \prod_{k=1}^n Q^p(\lambda_k)$$

$$+ \int_{B_n^c} d\lambda_1 \cdots d\lambda_n \prod_{k=1}^n [1 - \widehat{h}(\varepsilon(\lambda_k - \lambda_{k-1}))] Q^p(\lambda_k).$$

For the first term on the right-hand side,

$$\int_{B_n} d\lambda_1 \cdots d\lambda_n \prod_{k=1}^n Q^p(\lambda_k)$$

$$\leq \binom{n}{[n/3]} \left[ \int_{\mathbb{R}^d} Q^p(\lambda) \, d\lambda \right]^{n-[n/3]} \left[ \int_{\{|\lambda| \geq N\}} Q^p(\lambda) \, d\lambda \right]^{[n/3]}$$

$$\leq 2^n C^{n-[n/3]} \delta^{[n/3]}.$$

So we have

$$\limsup_{n \to \infty} \frac{1}{n} \log \int_{B_n} d\lambda_1 \cdots d\lambda_n \prod_{k=1}^n Q^p(\lambda_k) \leq \log 2 + \frac{2}{3} \log C + \frac{1}{3} \log \delta.$$

As for the second term, notice that on $B_n^c$, there are at least $[n/3]$ pairs $(\lambda_{k-1}, \lambda_k)$ such that $|\lambda_{k-1}| \leq N$ and $|\lambda_k| \leq N$. For such pairs we have $0 \leq 1 - \widehat{h}(\varepsilon(\lambda_k - \lambda_{k-1})) < \delta$. For any other pairs, we use the general bounds $0 \leq 1 - \widehat{h}(\varepsilon(\lambda_k - \lambda_{k-1})) \leq 1$. Therefore,

$$\int_{B_n^c} d\lambda_1 \cdots d\lambda_n \prod_{k=1}^n [1 - \widehat{h}(\varepsilon(\lambda_k - \lambda_{k-1}))] Q^p(\lambda_k)$$

$$\leq \delta^{[n/3]} \int_{(\mathbb{R}^d)^n} d\lambda_1 \cdots d\lambda_n \prod_{k=1}^n Q^p(\lambda_k) = C^n \delta^{[n/3]}.$$

Thus,

$$\limsup_{n\to\infty} \frac{1}{n} \log \int_{B_n^c} d\lambda_1 \cdots d\lambda_n \prod_{k=1}^{n} [1 - \widehat{h}(\varepsilon(\lambda_k - \lambda_{k-1}))] Q^p(\lambda_k)$$

$$\le \log C + \frac{1}{3} \log \delta.$$

In summary,

$$\limsup_{\varepsilon\to 0^+} \limsup_{n\to\infty} \frac{1}{n} \log \int_{(\mathbb{R}^d)^n} d\lambda_1 \cdots d\lambda_n \prod_{k=1}^{n} [1 - \widehat{h}(\varepsilon(\lambda_k - \lambda_{k-1}))] Q^p(\lambda_k)$$

$$\le \max\left\{ \log 2 + \frac{2}{3} \log C + \frac{1}{3} \log \delta, \ \log C + \frac{1}{3} \log \delta \right\}.$$

Letting $\delta \to 0^+$ gives (5.5).

By (5.4),

$$\limsup_{\varepsilon\to 0^+} \limsup_{n\to\infty} \frac{1}{n} \log \frac{1}{(n!)^p} \mathbb{E}[(\eta^0 - \eta_\varepsilon)([0,\tau_1] \times \cdots \times [0,\tau_p])^n] = -\infty.$$

We claim that it can be strengthened into

$$(5.6) \quad \limsup_{\varepsilon\to 0^+} \limsup_{n\to\infty} \frac{1}{n} \log \frac{1}{(n!)^p} \mathbb{E}|(\eta^0 - \eta_\varepsilon)([0,\tau_1] \times \cdots \times [0,\tau_p])|^n = -\infty.$$

Indeed, this is automatic if $n \to \infty$ along the even numbers. As for $n = 2k+1$, it is easy to see that our assertion follows from the following use of Hölder's inequality:

$$\mathbb{E}|(\eta^0 - \eta_\varepsilon)([0,\tau_1] \times \cdots \times [0,\tau_p])|^{2k+1}$$

$$\le \{\mathbb{E}|(\eta^0 - \eta_\varepsilon)([0,\tau_1] \times \cdots \times [0,\tau_p])|^{2k}\}^{1/2}$$

$$\times \{\mathbb{E}|(\eta^0 - \eta_\varepsilon)([0,\tau_1] \times \cdots \times [0,\tau_p])|^{2(k+1)}\}^{1/2}.$$

We now fix $\varepsilon > 0$ and estimate $\eta_\varepsilon([0,\tau_1] \times \cdots \times [0,\tau_p])$. Let $M > 0$ be fixed but arbitrary. By definition,

$$\eta_\varepsilon([0,t_1] \times \cdots \times [0,t_p])$$

$$(5.7) \qquad = \sum_{y\in\mathbb{Z}^d} \int_{[0,M]^d} h_\varepsilon(yM + z) \eta^{yM+z}([0,t_1] \times \cdots \times [0,t_p]) \, dz$$

$$\le \int_{[0,M]^d} \widetilde{h}_\varepsilon(z) \widetilde{\eta}^z([0,t_1] \times \cdots \times [0,t_p]) \, dz,$$

where

$$\widetilde{h}_\varepsilon(x) = \sum_{y\in\mathbb{Z}^d} h_\varepsilon(yM + z),$$

$$\widetilde{\eta}^z([0,t_1] \times \cdots \times [0,t_p]) = \sum_{y\in\mathbb{Z}^d} \eta^{yM+z}([0,t_1] \times \cdots \times [0,t_p])$$

are two periodic functions on $\mathbb{R}^d$ with the period $M > 0$.

By Parseval's identity,

$$\int_{[0,M]^d} \widetilde{h}_\varepsilon(z) \widetilde{\eta}^z([0,t_1] \times \cdots \times [0,t_p]) \, dz$$

$$= \frac{1}{M^d} \sum_{y \in \mathbb{Z}^d} \left( \int_{[0,M]^d} \widetilde{h}_\varepsilon(x) \exp\left\{ -i\frac{2\pi}{M}(y \cdot x) \right\} dx \right)$$

$$\times \left( \int_{[0,M]^d} \widetilde{\eta}^x([0,t_1] \times \cdots \times [0,t_p]) \exp\left\{ i\frac{2\pi}{M}(y \cdot x) \right\} dx \right).$$

By periodicity,

$$\int_{[0,M]^d} \widetilde{h}_\varepsilon(x) \exp\left\{ -i\frac{2\pi}{M}(y \cdot x) \right\} dx$$

$$= \sum_{z \in \mathbb{Z}^d} \int_{[0,M]^d} h_\varepsilon(zM + x) \exp\left\{ -i\frac{2\pi}{M}(y \cdot x) \right\} dx$$

$$= \sum_{z \in \mathbb{Z}^d} \int_{zM+[0,M]^d} h_\varepsilon(x) \exp\left\{ -i\frac{2\pi}{M}(y \cdot (x - zM)) \right\} dx$$

$$= \sum_{z \in \mathbb{Z}^d} \int_{zM+[0,M]^d} h_\varepsilon(x) \exp\left\{ -i\frac{2\pi}{M}(y \cdot x) \right\} dx$$

$$= \int_{\mathbb{R}^d} h_\varepsilon(x) \exp\left\{ -i\frac{2\pi}{M}(y \cdot x) \right\} dx = \widehat{h}\left( \frac{2\pi\varepsilon}{M} y \right).$$

Similarly,

$$\int_{[0,M]^d} \widetilde{\eta}^x([0,t_1] \times \cdots \times [0,t_p]) \exp\left\{ i\frac{2\pi}{M}(y \cdot x) \right\} dx$$

$$= \int_{\mathbb{R}^d} \eta^x([0,t_1] \times \cdots \times [0,t_p]) \exp\left\{ i\frac{2\pi}{M}(y \cdot x) \right\} dx$$

$$= \int_{[0,t_1] \times \cdots \times [0,t_p]} \exp\left\{ i\frac{2\pi}{M} y \cdot (X_1(s_1) + \cdots + X_p(s_p)) \right\} ds_1 \cdots ds_p.$$

Hence,

$$\int_{[0,M]^d} \widetilde{h}_\varepsilon(z) \widetilde{\eta}^z([0,t_1] \times \cdots \times [0,t_p]) \, dz$$

$$= \frac{1}{M^d} \sum_{y \in \mathbb{Z}^d} \widehat{h}\left( \frac{2\pi\varepsilon}{M} y \right)$$

(5.8)

$$\times \int_{[0,t_1] \times \cdots \times [0,t_p]} \exp\left\{ i\frac{2\pi}{M} y \right.$$

$$\cdot\, (X_1(s_1) + \cdots + X_p(s_p))\Big\}\, ds_1 \cdots ds_p.$$

Following a procedure same as the one used for (2.2),

$$\mathbb{E}\left[\int_{[0,M]^d} \widetilde{h}_\varepsilon(z)\widetilde{\eta}^z([0,\tau_1] \times \cdots \times [0,\tau_p])\, dz\right]^n$$

$$= \frac{1}{M^{dn}} \sum_{y_1,\ldots,y_n \in \mathbb{Z}^d} \left(\prod_{k=1}^n \widehat{h}\left(\frac{2\pi\varepsilon}{M}y_k\right)\right)\left[\sum_{\sigma \in \Sigma_n} \prod_{k=1}^n Q\left(\frac{2\pi}{M}\sum_{j=1}^k y_{\sigma(j)}\right)\right]^p.$$

By Theorem 4.1,

$$\lim_{n\to\infty} \frac{1}{n}\log\frac{1}{(n!)^p}\mathbb{E}\left[\int_{[0,M]^d} \widetilde{h}_\varepsilon(z)\widetilde{\eta}^z([0,\tau_1] \times \cdots \times [0,\tau_p])\, dz\right]^n$$

$$= \log\left(\frac{1}{M^d} \sup_{|f|_2=1} \sum_{x\in\mathbb{Z}^d} \widehat{h}\left(\frac{2\pi\varepsilon}{M}x\right)\left[\sum_{y\in\mathbb{Z}^d} \sqrt{Q\left(\frac{2\pi}{M}(x+y)\right)}\right.\right.$$

(5.9)
$$\left.\left. \times \sqrt{Q\left(\frac{2\pi}{M}y\right)}f(x+y)f(y)\right]^p\right)$$

$$\leq \log(M^{-d}\rho_M),$$

where

$$(5.10) \quad \rho_M = \sup_{|f|_2=1} \sum_{x\in\mathbb{Z}^d}\left[\sum_{y\in\mathbb{Z}^d} \sqrt{Q\left(\frac{2\pi}{M}(x+y)\right)}\sqrt{Q\left(\frac{2\pi}{M}y\right)}f(x+y)f(y)\right]^p.$$

In view of (5.7),

$$\limsup_{n\to\infty} \frac{1}{n}\log\frac{1}{(n!)^p}\mathbb{E}[\eta_\varepsilon([0,\tau_1] \times \cdots \times [0,\tau_p])]^n \leq \log(M^{-d}\rho_M).$$

By Lemma A.1 given in the Appendix, letting $M \to \infty$ on the right-hand side gives

$$(5.11) \quad \limsup_{n\to\infty} \frac{1}{n}\log\frac{1}{(n!)^p}\mathbb{E}[\eta_\varepsilon([0,\tau_1] \times \cdots \times [0,\tau_p])]^n \leq \log\frac{\rho}{(2\pi)^d}.$$

Finally, (5.1) follows from (5.6), (5.11) and the fact that

$$\{\mathbb{E}[\eta^0([0,\tau_1] \times \cdots \times [0,\tau_p])]^n\}^{1/n}$$

$$\leq \{\mathbb{E}[\eta_\varepsilon([0,\tau_1] \times \cdots \times [0,\tau_p])]^n\}^{1/n}$$

$$+ \{\mathbb{E}|(\eta^0 - \eta_\varepsilon)([0,\tau_1] \times \cdots \times [0,\tau_p])|^n\}^{1/n}.$$

**6. Proof of Theorem 1.2.** In the light of Theorem 1.1, the nontrivial part of Theorem 1.2 is the upper bound. Let $M > 0$ be fixed and recall that

$$\widetilde{\eta}^x([0, \tau_1] \times \cdots \times [0, \tau_p]) = \sum_{y \in \mathbb{Z}^d} \eta^{yM+x}([0, \tau_1] \times \cdots \times [0, \tau_p]).$$

Notice that

$$(6.1) \quad \sup_{x \in \mathbb{R}^d} \eta^x([0, \tau_1] \times \cdots \times [0, \tau_p]) \leq \sup_{x \in [0, M]^d} \widetilde{\eta}^x([0, \tau_1] \times \cdots \times [0, \tau_p]).$$

By Fourier expansion,

$$\widetilde{\eta}^x([0, \tau_1] \times \cdots \times [0, \tau_p]) = \sum_{y \in \mathbb{Z}^d} a(y) \exp\left\{ i \frac{2\pi}{M}(x \cdot y) \right\},$$

where

$$
\begin{aligned}
a(y) &= \frac{1}{M^d} \int_{[0, M]^d} \exp\left\{ -i \frac{2\pi}{M}(x \cdot y) \right\} \widetilde{\eta}^x([0, \tau_1] \times \cdots \times [0, \tau_p]) \, dx \\
&= \frac{1}{M^d} \int_{\mathbb{R}^d} \exp\left\{ -i \frac{2\pi}{M}(x \cdot y) \right\} \eta^x([0, \tau_1] \times \cdots \times [0, \tau_p]) \, dx \\
&= \frac{1}{M^d} \int_0^{\tau_1} \cdots \int_0^{\tau_p} \exp\left\{ -i \left( \frac{2\pi}{M} y \right) \cdot (X_1(s_1) + \cdots + X_p(s_p)) \right\} ds_1 \cdots ds_p.
\end{aligned}
$$

Thus

$$
\begin{aligned}
&\widetilde{\eta}^x([0, \tau_1] \times \cdots \times [0, \tau_p]) \\
(6.2) \quad &= \frac{1}{M^d} \sum_{y \in \mathbb{Z}^d} \exp\left\{ i \frac{2\pi}{M}(x \cdot y) \right\} \\
&\qquad\qquad \times \int_0^{\tau_1} \cdots \int_0^{\tau_p} \exp\left\{ -i \left( \frac{2\pi}{M} y \right) \right. \\
&\qquad\qquad\qquad\qquad \left. \cdot (X_1(s_1) + \cdots + X_p(s_p)) \right\} ds_1 \cdots ds_p.
\end{aligned}
$$

Let the functions $h$ and $h_\varepsilon$ be defined in (5.2) and (5.3), respectively. Recall that

$$\widetilde{h}_\varepsilon(x) = \sum_{y \in \mathbb{Z}^d} h_\varepsilon(yM + x), \qquad \widehat{h}(\lambda) = \int_{\mathbb{R}^d} h(x) e^{i\lambda \cdot x} \, dx.$$

Write

$$\widetilde{\eta}_\varepsilon([0, \tau_1] \times \cdots \times [0, \tau_p]) = \int_{[0, M]^d} \widetilde{h}_\varepsilon(x) \widetilde{\eta}^x([0, \tau_1] \times \cdots \times [0, \tau_p]) \, dx.$$

By (5.8) and (6.2), and by a procedure similar to the one for (2.2), one can prove that

$$\mathbb{E}[\widetilde{\eta}^x([0,\tau_1] \times \cdots \times [0,\tau_p])]^n$$

$$(6.3) \qquad = \frac{1}{M^{dn}} \sum_{y_1,\ldots,y_n \in \mathbb{Z}^d} \exp\left\{ i \frac{2\pi}{M} \sum_{k=1}^{n} (x \cdot y_k) \right\}$$

$$\times \left[ \sum_{\sigma \in \Sigma_n} \prod_{k=1}^{n} Q\left( \frac{2\pi}{M} \sum_{j=1}^{k} y_{\sigma(j)} \right) \right]^p,$$

$$\mathbb{E}[(\widetilde{\eta}^x - \widetilde{\eta}^z)([0,\tau_1] \times \cdots \times [0,\tau_p])]^n$$

$$(6.4) \qquad = \frac{1}{M^{dn}} \sum_{y_1,\ldots,y_n \in \mathbb{Z}^d} \left( \prod_{k=1}^{n} \left[ \exp\left\{ i \frac{2\pi}{M} (x \cdot y_k) \right\} - \exp\left\{ i \frac{2\pi}{M} (z \cdot y_k) \right\} \right] \right)$$

$$\times \left[ \sum_{\sigma \in \Sigma_n} \prod_{k=1}^{n} Q\left( \frac{2\pi}{M} \sum_{j=1}^{k} y_{\sigma(j)} \right) \right]^p,$$

$$x, z \in [0, M]^d,$$

$$\mathbb{E}[(\widetilde{\eta}^0 - \widetilde{\eta}_\varepsilon)([0,\tau_1] \times \cdots \times [0,\tau_p])]^n$$

$$(6.5) \qquad = \frac{1}{M^{dn}} \sum_{y_1,\ldots,y_n \in \mathbb{Z}^d} \left( \prod_{k=1}^{n} \left[ 1 - \widehat{h}\left( \frac{2\pi\varepsilon}{M} y_k \right) \right] \right)$$

$$\times \left[ \sum_{\sigma \in \Sigma_n} \prod_{k=1}^{n} Q\left( \frac{2\pi}{M} \sum_{j=1}^{k} y_{\sigma(j)} \right) \right]^p.$$

By (6.5) and by an argument similar to the one used for (5.6),

$$\limsup_{\varepsilon \to 0^+} \limsup_{n \to \infty} \frac{1}{n} \log \frac{1}{(n!)^p} \mathbb{E}|(\widetilde{\eta}^0 - \widetilde{\eta}_\varepsilon)([0,\tau_1] \times \cdots \times [0,\tau_p])|^n = -\infty.$$

This, together with (5.9), implies that

$$(6.6) \qquad \limsup_{n \to \infty} \frac{1}{n} \log \frac{1}{(n!)^p} \mathbb{E}[\widetilde{\eta}^0([0,\tau_1] \times \cdots \times [0,\tau_p])]^n \le \log \frac{\rho_M}{M^d}.$$

By Lemma 6.1 given below and by Taylor's expansion one can easily see that

$$(6.7) \qquad \begin{aligned} &\limsup_{\delta \to 0^+} \limsup_{n \to \infty} \frac{1}{n} \log \frac{1}{(n!)^p} \mathbb{E} \sup_{|y-x| \le \delta} |(\widetilde{\eta}^y - \widetilde{\eta}^x)([0,\tau_1] \times \cdots \times [0,\tau_p])|^n \\ &= -\infty. \end{aligned}$$

Given $\delta > 0$, let $D \subset [0,M]^d$ be a finite $\delta$-net of $[0,M]^d$:

$$\left\{ \mathbb{E} \sup_{x \in [0,M]^d} \widetilde{\eta}^x([0,\tau_1] \times \cdots \times [0,\tau_p])^n \right\}^{1/n}$$

$$\leq \left\{ \mathbb{E} \sup_{x \in D} \widetilde{\eta}^x([0,\tau_1] \times \cdots \times [0,\tau_p])^n \right\}^{1/n}$$

(6.8)
$$+ \left\{ \mathbb{E} \sup_{|y-x| \leq \delta} |(\widetilde{\eta}^y - \widetilde{\eta}^x)([0,\tau_1] \times \cdots \times [0,\tau_p])|^n \right\}^{1/n}$$

$$\leq \left\{ \#(D) \sup_{x \in [0,M]^d} \mathbb{E} \widetilde{\eta}^x([0,\tau_1] \times \cdots \times [0,\tau_p])^n \right\}^{1/n}$$

$$+ \left\{ \mathbb{E} \sup_{|y-x| \leq \delta} |(\widetilde{\eta}^y - \widetilde{\eta}^x)([0,\tau_1] \times \cdots \times [0,\tau_p])|^n \right\}^{1/n}.$$

From (6.3) one can see that for any $x \in [0,M]^d$ and for any integer $n \geq 0$,

$$\mathbb{E}[\widetilde{\eta}^x([0,\tau_1] \times \cdots \times [0,\tau_p])^n] \leq \mathbb{E}[\widetilde{\eta}^0([0,\tau_1] \times \cdots \times [0,\tau_p])^n].$$

By (6.6), (6.7) and (6.8), therefore,

$$\limsup_{n \to \infty} \frac{1}{n} \log \frac{1}{(n!)^p} \mathbb{E} \left[ \sup_{x \in [0,M]^d} \widetilde{\eta}^x([0,\tau_1] \times \cdots \times [0,\tau_p])^n \right] \leq \log \frac{\rho_M}{M^d}.$$

In view of (6.1), we have

$$\limsup_{n \to \infty} \frac{1}{n} \log \frac{1}{(n!)^p} \mathbb{E} \left[ \sup_{x \in \mathbb{R}^d} \eta^x([0,\tau_1] \times \cdots \times [0,\tau_p])^n \right] \leq \log \frac{\rho_M}{M^d}.$$

By Lemma A.1 given in the Appendix below, letting $M \to \infty$ on the right-hand side gives

(6.9)   $$\limsup_{n \to \infty} \frac{1}{n} \log \frac{1}{(n!)^p} \mathbb{E} \left[ \sup_{x \in \mathbb{R}^d} \eta^x([0,\tau_1] \times \cdots \times [0,\tau_p])^n \right] \leq \log \frac{\rho}{(2\pi)^d}.$$

We now adopt the argument used for (2.5) here. For this we replace (2.3) by (6.9), and replace (1.6) by (1.13). We obtain

$$\limsup_{n \to \infty} \frac{1}{n} \log(n!)^{-d/\alpha} \mathbb{E} \left[ \sup_{x \in \mathbb{R}^d} \eta^x([0,1]^p)^n \right]$$

$$\leq \log \left( \frac{\alpha p}{\alpha p - d} \right)^{(\alpha p - d)/\alpha} + \log \frac{\rho}{(2\pi)^d}.$$

Comparing this to (2.6) gives

$$\lim_{n \to \infty} \frac{1}{n} \log(n!)^{-d/\alpha} \mathbb{E} \left[ \sup_{x \in \mathbb{R}^d} \eta^x([0,1]^p)^n \right]$$

(6.10)
$$= \log\left(\frac{\alpha p}{\alpha p - d}\right)^{(\alpha p - d)/\alpha} + \log\frac{\rho}{(2\pi)^d}.$$

Finally, Theorem 1.2 follows from (6.10) and Lemma 1.4.

LEMMA 6.1.   *For any number $\zeta$ with $0 < \zeta < \min\{1, \ (\alpha p - d)/2\}$, there is a positive number $c = c(\zeta, \psi, p)$ such that*

(6.11) $\mathbb{E}\exp\left\{ c \sup_{\substack{x,z\in[0,M]^d \\ x\neq z}} \left(\frac{|(\widetilde{\eta}^x - \widetilde{\eta}^z)([0,\tau_1]\times\cdots\times[0,\tau_p])|}{|x-z|^\zeta}\right)^{1/p}\right\} < \infty.$

PROOF.   By (6.4) and Jensen's inequality, for any $x, z \in [0,M]^d$,

$\mathbb{E}[(\widetilde{\eta}^x - \widetilde{\eta}^z)([0,\tau_1]\times\cdots\times[0,\tau_p])^n]$

$$\leq \frac{1}{M^{dn}} \sum_{y_1,\ldots,y_k\in\mathbb{Z}^d} \left(\prod_{k=1}^n \left|1 - \exp\left\{i\frac{2\pi}{M}((z-x)\cdot y_k)\right\}\right|\right)$$
$$\times \left[\sum_{\sigma\in\Sigma_n}\prod_{k=1}^n Q\left(\frac{2\pi}{M}\sum_{j=1}^k y_{\sigma(j)}\right)\right]^p$$

$$\leq \frac{(n!)^p}{M^{dn}} \sum_{y_1,\ldots,y_k\in\mathbb{Z}^d} \left(\prod_{k=1}^n \left|1 - \exp\left\{i\frac{2\pi}{M}((z-x)\cdot y_k)\right\}\right|\right) \prod_{k=1}^n Q^p\left(\frac{2\pi}{M}\sum_{j=1}^k y_j\right)$$

$$= \frac{(n!)^p}{M^{dn}} \sum_{y_1,\ldots,y_k\in\mathbb{Z}^d} \left(\prod_{k=1}^n \left|1 - \exp\left\{i\frac{2\pi}{M}((z-x)\cdot(y_k - y_{k-1}))\right\}\right|\right)$$
$$\times \prod_{k=1}^n Q^p\left(\frac{2\pi}{M}y_k\right).$$

Fix $\zeta'$ with $\zeta < \zeta' < \min\{1, \ (\alpha p - d)/2\}$. We have

(6.12) $$\int_{\mathbb{R}^d} |\lambda|^\theta Q^p(\lambda)\,d\lambda < \infty, \qquad \theta \leq 2\zeta'.$$

Notice that

$$\left|1 - \exp\left\{i\frac{2\pi}{M}((z-x)\cdot(y_k - y_{k-1}))\right\}\right| \leq 1.$$

Hence,

$$\frac{1}{2}\left|1 - \exp\left\{i\frac{2\pi}{M}((z-x)\cdot(y_k - y_{k-1}))\right\}\right|$$

$$\leq 2^{-\zeta'} \left| 1 - \exp\left\{ i\frac{2\pi}{M}((z-x)\cdot(y_k - y_{k-1})) \right\} \right|^{\zeta'}$$

$$\leq 2^{-\zeta'} \left(\frac{2\pi}{M}\right)^{\zeta'} |z-x|^{\zeta'} |y_k - y_{k-1}|^{\zeta'}.$$

Therefore,

$$
\begin{aligned}
(6.13) \quad & \mathbb{E}[(\widetilde{\eta}^x - \widetilde{\eta}^z)([0,\tau_1] \times \cdots \times [0,\tau_p])^n] \\
& \leq 2^n \frac{(n!)^p}{M^{nd}} \left(\frac{2\pi}{M}|x-z|\right)^{\zeta' n} \sum_{y_1,\ldots,y_k \in \mathbb{Z}^d} \prod_{k=1}^{n} |y_k - y_{k-1}|^{\zeta'} Q^p\left(\frac{2\pi}{M}y_k\right).
\end{aligned}
$$

By the triangular inequality,

$$\prod_{k=1}^{n} |y_k - y_{k-1}|^{\zeta'} \leq \prod_{k=1}^{n} (|y_k|^{\zeta'} + |y_{k-1}|^{\zeta'}) = \sum_{\delta_1,\ldots,\delta_n} \prod_{k=1}^{n} |y_k|^{\delta_k \zeta'},$$

where for each $1 \leq k \leq n$, $\delta_k$ has three possible values: 0, 1, or 2, and $\delta_1 + \cdots + \delta_n = n$. The total number of the terms is at most $2^n$. Thus,

$$
\begin{aligned}
& \left(\frac{2\pi}{M}\right)^{nd} \left(\frac{2\pi}{M}\right)^{\zeta' n} \sum_{y_1,\ldots,y_k \in \mathbb{Z}^d} \prod_{k=1}^{n} |y_k - y_{k-1}|^{\zeta'} Q^p\left(\frac{2\pi}{M}y_k\right) \\
& \leq \sum_{\delta_1,\ldots,\delta_n} \prod_{k=1}^{n} \left(\frac{2\pi}{M}\right)^d \sum_{y \in \mathbb{Z}^d} \left|\frac{2\pi}{M}y\right|^{\delta_k \zeta'} Q^p\left(\frac{2\pi}{M}y\right).
\end{aligned}
$$

From (6.12) there is a $C = C(\zeta, \psi, p) > 0$ such that

$$\left(\frac{2\pi}{M}\right)^d \sum_{y \in \mathbb{Z}^d} \left|\frac{2\pi}{M}y\right|^{\delta_k \zeta'} Q^p\left(\frac{2\pi}{M}y\right) \leq C.$$

So we have

$$
\begin{aligned}
& \left(\frac{2\pi}{M}\right)^{nd} \left(\frac{2\pi}{M}\right)^{\zeta' n} \sum_{y_1,\ldots,y_k \in \mathbb{Z}^d} \prod_{k=1}^{n} |y_k - y_{k-1}|^{\zeta'} Q^p\left(\frac{2\pi}{M}y_k\right) \\
& \leq \sum_{\delta_1,\ldots,\delta_n} \prod_{k=1}^{n} C \leq (2C)^n.
\end{aligned}
$$

By (6.13),

$$\sup_{\substack{x,z \in [0,M]^d \\ x \neq z}} \mathbb{E}\left[\frac{(\widetilde{\eta}^x - \widetilde{\eta}^z)([0,\tau_1] \times \cdots \times [0,\tau_p])}{|x-z|^{\zeta'}}\right]^n \leq (n!)^p (4C)^n,$$

$$n = 0, 1, 2, \ldots.$$

Following a standard way of using Hölder inequality, we conclude that there is a $C_0 = C_0(\zeta, \psi, p) > 0$ such that

$$\sup_{\substack{x, z \in [0, M]^d \\ x \neq z}} \mathbb{E} \left| \frac{(\widetilde{\eta}^x - \widetilde{\eta}^z)([0, \tau_1] \times \cdots \times [0, \tau_p])}{|x - z|^{\zeta'}} \right|^n \leq (n!)^p C_0^n,$$

(6.14)

$$n = 0, 1, 2, \ldots.$$

Recall that a function $\Psi : \mathbb{R}^+ \longrightarrow \mathbb{R}^+$ is a Young function if it is convex, increasing and satisfies $\Psi(0) = 0$, $\lim_{x \to \infty} \Psi(x) = \infty$. The Orlicz space $\mathcal{L}_\Psi(\Omega, \mathcal{A}, \mathbb{P})$ is defined as the linear space of all random variables $X$ on the probability space $(\Omega, \mathcal{A}, \mathbb{P})$ such that

$$\|X\|_\Psi = \inf\{c > 0; \ \mathbb{E}\Psi(c^{-1}|X|) \leq 1\} < \infty.$$

It has been known that $\|\cdot\|_\Psi$ defines a norm (called Orlicz norm) and $\mathcal{L}_\Psi(\Omega, \mathcal{A}, \mathbb{P})$ becomes a Banach space under $\|\cdot\|_\Psi$.

We now choose the Young function $\Psi$ such that $\Psi(x) \sim \exp\{x^{1/p}\}$ as $x \to \infty$. By (6.14) there is $c = c(\zeta, d, p) > 0$ such that

$$\|(\widetilde{\eta}^x - \widetilde{\eta}^z)([0, \tau_1] \times \cdots \times [0, \tau_p])\|_\Psi \leq c|x - z|^{\zeta'}, \qquad x, z \in [0, M]^d.$$

By a standard chaining argument (see, e.g., Lemma 9 in [7]),

$$\left\| \sup_{\substack{x, z \in [0, M]^d \\ x \neq z}} \frac{|(\widetilde{\eta}^x - \widetilde{\eta}^z)([0, \tau_1] \times \cdots \times [0, \tau_p])|}{|x - z|^\zeta} \right\|_\Psi < \infty,$$

which leads to the desired conclusion. $\square$

**7. Proof of Theorem 1.3.** The upper bound in (1.16) and therefore the upper bound in (1.15) follow from Theorem 1.2, the scaling property given in (1.13) and a standard procedure via the Borel–Cantelli lemma. It remains to prove that for any fix $x \in \mathbb{R}^d$,

$$\limsup_{t \to \infty} t^{-(\alpha p - d)/\alpha} (\log \log t)^{-d/\alpha} \eta^x([0, t]^p)$$

(7.1)

$$\geq (2\pi)^{-d} \left( \frac{\alpha}{d} \right)^{d/\alpha} \left( 1 - \frac{d}{\alpha p} \right)^{-(p - d/\alpha)} \rho \qquad \text{a.s.}$$

We first prove that

$$\lim_{\delta \to 0^+} \liminf_{t \to \infty} t^{-1} \log \mathbb{P} \left\{ \inf_{|y| \leq \delta} \eta^y([0, t]^p) \geq t^p \right\}$$

(7.2)

$$\geq -(2\pi)^\alpha \frac{d}{\alpha} \left( 1 - \frac{d}{\alpha p} \right)^{(\alpha p - d)/d} \rho^{-\alpha/d}.$$

Indeed, similarly to Lemma 6.1, for any bounded neighborhood $D$ of 0 and any $0 < \zeta < \min\{1, (\alpha p - d)/2\}$ there is a $c = c(D, \zeta, \psi, p) > 0$ such that

$$(7.3) \qquad \mathbb{E} \exp\left\{ c \sup_{\substack{y, z \in D \\ y \neq z}} \left( \frac{|(\eta^y - \eta^z)([0, \tau_1] \times \cdots \times [0, \tau_p])|}{|y - z|^{\zeta}} \right)^{1/p} \right\} < \infty.$$

By the Chebyshev inequality we have that for any $\varepsilon > 0$,

$$\limsup_{\delta \to 0^+} \limsup_{t \to \infty} t^{-1} \log \mathbb{P}\left\{ \sup_{|y| \leq \delta} |(\eta^0 - \eta^y)([0, \tau_1] \times \cdots \times [0, \tau_p])| \geq \varepsilon t^p \right\}$$
$$= -\infty.$$

On the other hand,

$$\mathbb{P}\left\{ \sup_{|y| \leq \delta} |(\eta^0 - \eta^y)([0, \tau_1] \times \cdots \times [0, \tau_p])| \geq \varepsilon t^p \right\}$$
$$= \int_0^\infty \cdots \int_0^\infty e^{-(t_1 + \cdots + t_p)}$$
$$\times \mathbb{P}\left\{ \sup_{|y| \leq \delta} |(\eta^0 - \eta^y)([0, t_1] \times \cdots \times [0, t_p])| \geq \varepsilon t^p \right\} dt_1 \cdots dt_p$$
$$\geq \int_{(1-\varepsilon)t}^t \cdots \int_{(1-\varepsilon)t}^t e^{-(t_1 + \cdots + t_p)}$$
$$\times \mathbb{P}\left\{ \sup_{|y| \leq \delta} |(\eta^0 - \eta^y)([0, t_1]\right.$$
$$\left. \times \cdots \times [0, t_p])| \geq \varepsilon t^p \right\} dt_1 \cdots dt_p$$
$$\geq (e^{-(1-\varepsilon)t} - e^{-t})^p$$
$$\times \inf_{(1-\varepsilon)t \leq t_1, \ldots, t_p \leq t} \mathbb{P}\left\{ \sup_{|y| \leq \delta} |(\eta^0 - \eta^y)([0, t_1] \times \cdots \times [0, t_p])| \geq \varepsilon t^p \right\}.$$

So we have

$$(7.4) \qquad \begin{aligned} &\limsup_{\delta \to 0^+} \limsup_{t \to \infty} t^{-1} \log \inf_{(1-\varepsilon)t \leq t_1, \ldots, t_p \leq t} \\ &\mathbb{P}\left\{ \sup_{|y| \leq \delta} |(\eta^0 - \eta^y)([0, t_1] \times \cdots \times [0, t_p])| \geq \varepsilon t^p \right\} = -\infty. \end{aligned}$$

For any $t$ and $(1-\varepsilon)t \leq t_1, \ldots, t_p \leq t$,

$$\inf_{|y| \leq \delta} \eta^y([0, t]^p)$$

$$\geq \inf_{|y|\leq\delta} \eta^y([0,t_1]\times\cdots\times[0,t_p])$$

$$\geq \eta^0([0,t_1]\times\cdots\times[0,t_p]) - \inf_{|y|\leq\delta}|(\eta^0-\eta^y)([0,t_1]\times\cdots\times[0,t_p])|$$

$$\geq \eta^0([0,(1-\varepsilon)t]^p) - \inf_{|y|\leq\delta}|(\eta^0-\eta^y)([0,t_1]\times\cdots\times[0,t_p])|.$$

Hence,

$$\mathbb{P}\Big\{\inf_{|x|\leq\delta}\eta^x([0,t]^p)\geq t^p\Big\}$$

$$+ \inf_{(1-\varepsilon)t\leq t_1,\ldots,t_p\leq t}\mathbb{P}\Big\{\sup_{|x|\leq\delta}|(\eta^0-\eta^x)([0,t_1]\times\cdots\times[0,t_p])|\geq\varepsilon t^p\Big\}$$

$$\geq \mathbb{P}\{\eta^0([0,(1-\varepsilon)t]^p)\geq(1+\varepsilon)t^p\}.$$

Consequently,

$$(7.5) \quad \begin{aligned} &\max\Big\{\liminf_{t\to\infty}t^{-1}\log\mathbb{P}\Big\{\inf_{|y|\leq\delta}\eta^y([0,t]^p)\geq t^p\Big\},\\ &\qquad\limsup_{t\to\infty}t^{-1/p}\log\inf_{(1-\varepsilon)t\leq t_1,\ldots,t_p\leq t}\\ &\qquad\mathbb{P}\Big\{\sup_{|y|\leq\delta}|(\eta^0-\eta^y)([0,t_1]\times\cdots\times[0,t_p])|\geq\varepsilon t\Big\}\Big\}\\ &\qquad\geq\lim_{t\to\infty}t^{-1}\log\mathbb{P}\{\eta^0([0,(1-\varepsilon)t]^p)\geq(1+\varepsilon)t^p\}. \end{aligned}$$

Notice that

$$\mathbb{P}\{\eta^0([0,(1-\varepsilon)t]^p)\geq(1+\varepsilon)t^p\}$$

$$= \mathbb{P}\{\eta^0([0,1]^p)\geq(1+\varepsilon)(1-\varepsilon)^{-(\alpha p-d)/\alpha}t^{d/\alpha}\}.$$

By Theorem 1.1,

$$(7.6) \quad \begin{aligned} &\lim_{t\to\infty}t^{-1}\log\mathbb{P}\{\eta^0([0,(1-\varepsilon)t]^p)\geq(1-\varepsilon)t^p\}\\ &= -(1+\varepsilon)^{\alpha/d}(1-\varepsilon)^{-(\alpha p-d)-d}(2\pi)^\alpha\frac{d}{\alpha}\Big(1-\frac{d}{\alpha p}\Big)^{(\alpha p-d)/d}\rho^{-\alpha/d}. \end{aligned}$$

Let $\delta\to0^+$ in (7.5). By (7.4) and (7.6) we obtain

$$\lim_{\delta\to0^+}\liminf_{t\to\infty}t^{-1}\log\mathbb{P}\Big\{\inf_{|y|\leq\delta}\eta^y([0,t]^p)\geq t^p\Big\}$$

$$\geq -(1+\varepsilon)^{\alpha/d}(1-\varepsilon)^{-(\alpha p-d)/d}(2\pi)^\alpha\frac{d}{\alpha}\Big(1-\frac{d}{\alpha p}\Big)^{(\alpha p-d)/d}\rho^{-\alpha/d}.$$

Letting $\varepsilon \to 0^+$ on the right-hand side leads to (7.2).

We come to the proof of (7.1). For each $k \geq 1$, write $t_k = k^k$ and define

$$X_{j,k}(t) = X_j(t_k + t) - X_j(t_k), \qquad t \geq 0, \; j = 1, \ldots, p, \; k = 1, 2, \ldots.$$

Let $\eta_k^x(I)$ be the local time of the additive stable process

$$\overline{X}_k(s_1, \ldots, s_p) = X_{1,k}(s_1) + \cdots + X_{p,k}(s_p).$$

Then for each $k$, $\eta_k \stackrel{d}{=} \eta$.

Let $\delta > 0$ be a small number which will be specified later. Write $Y_k = X_1(t_k) + \cdots + X_p(t_k)$. A rough estimate gives that with probability 1, the inequality

$$|Y_k| \leq 2^{-1} \delta \left( \frac{t_{k+1}}{\log t_{k+1}} \right)^{1/\alpha}$$

eventually holds. Therefore, with probability 1,

$$(7.7) \quad \begin{aligned} \eta^x([t_k, t_{k+1}]^p) &= \eta_k^{x+Y_k}([0, t_{k+1} - t_k]^p) \\ &\geq \inf_{|y| \leq \delta(t_{k+1}/\log\log t_{k+1})^{1/\alpha}} \eta_k^y([0, t_{k+1} - t_k]^p) \end{aligned}$$

eventually holds.

For each $k$, by the scaling property of the stable processes,

$$\begin{aligned} &\inf_{|y| \leq \delta(t_{k+1}/\log\log t_{k+1})^{1/\alpha}} \eta_k^y([0, t_{k+1} - t_k]^p) \\ &\stackrel{d}{=} \inf_{|y| \leq \delta(t_{k+1}/\log\log t_{k+1})^{1/\alpha}} \eta^y([0, t_{k+1} - t_k]^p) \\ &\stackrel{d}{=} \left( \frac{t_{k+1}}{\log\log t_{k+1}} \right)^{(\alpha p - d)/\alpha} \inf_{|y| \leq \delta} \eta^y([0, \; t_{k+1}^{-1}(t_{k+1} - t_k) \log\log t_{k+1}]^p). \end{aligned}$$

Let $\theta > 0$ satisfy

$$\theta < (2\pi)^{-d} \left( \frac{\alpha}{d} \right)^{d/\alpha} \left( 1 - \frac{d}{\alpha p} \right)^{-(p - d/\alpha)} \rho.$$

We have

$$\begin{aligned} &\mathbb{P}\left\{ \inf_{|y| \leq \delta(t_{k+1}/\log\log t_{k+1})^{1/\alpha}} \eta_k^y([0, t_{k+1} - t_k]^p) \geq \theta t_{k+1}^{(\alpha p - d)/\alpha} (\log\log t_{k+1})^{d/\alpha} \right\} \\ &= \mathbb{P}\left\{ \inf_{|y| \leq \delta} \eta^y([0, \; t_{k+1}^{-1}(t_{k+1} - t_k) \log\log t_{k+1}]^p) \geq \theta(\log\log t_{k+1})^p \right\}. \end{aligned}$$

By (7.2), therefore, one can take $\delta > 0$ sufficiently small so that

$$\liminf_{k \to \infty} \frac{1}{\log \log t_{k+1}} \log \mathbb{P} \left\{ \inf_{|y| \leq \delta(t_{k+1}/\log \log t_{k+1})^{1/\alpha}} \eta_k^y([0, t_{k+1} - t_k]^p) \right.$$

$$\left. \geq \theta t_{k+1}^{(\alpha p - d)/\alpha} (\log \log t_{k+1})^{d/\alpha} \right\} > -1.$$

Consequently,

$$\sum_k \mathbb{P} \left\{ \inf_{|y| \leq \delta(t_{k+1}/\log \log t_{k+1})^{1/\alpha}} \eta_k^y([0, t_{k+1} - t_k]^p) \right.$$

$$\left. \geq \theta t_{k+1}^{(\alpha p - d)/\alpha} (\log \log t_{k+1})^{d/\alpha} \right\} = \infty.$$

Notice that

$$\inf_{|y| \leq \delta(t_{k+1}/\log \log t_{k+1})^{1/\alpha}} \eta_k^y([0, t_{k+1} - t_k]^p), \qquad k = 1, 2, \ldots$$

is an independent sequence. By the Borel–Cantelli lemma,

$$\limsup_{k \to \infty} t_{k+1}^{-(\alpha p - d)/\alpha} (\log \log t_{k+1})^{-d/\alpha}$$

$$\times \inf_{|y| \leq \delta(t_{k+1}/\log \log t_{k+1})^{1/\alpha}} \eta_k^y([0, t_{k+1} - t_k]^p) \geq \theta \qquad \text{a.s.}$$

By (7.7),

$$\limsup_{k \to \infty} t_{k+1}^{-(\alpha p - d)/\alpha} (\log \log t_{k+1})^{-d/\alpha} \eta^x([t_k, t_{k+1}]^p) \geq \theta \qquad \text{a.s.}$$

Consequently,

$$\limsup_{t \to \infty} t^{-(\alpha p - d)/\alpha} (\log \log t)^{-d/\alpha} \eta^x([0, t]^p) \geq \theta \qquad \text{a.s.}$$

Letting

$$\theta \longrightarrow (2\pi)^{-d} \left( \frac{\alpha}{d} \right)^{d/\alpha} \left( 1 - \frac{d}{\alpha p} \right)^{-(p - d/\alpha)} \rho^-$$

proves (7.1).

## APPENDIX

LEMMA A.1. *Let $\rho$ be defined in (1.4) and let $\rho_M$ be defined in (5.10). We have*

$$\text{(A.1)} \qquad \limsup_{M \to \infty} M^{-d} \rho_M \leq (2\pi)^{-d} \rho.$$

PROOF.    Given $a > 0$ and $f \in \mathcal{L}^2(\mathbb{Z}^d)$ with $|f|_2 = 1$, by Hölder inequality

$$\sum_{|x| \geq (2\pi)^{-1}Ma} \left[ \sum_{y \in \mathbb{Z}^d} \sqrt{Q\left(\frac{2\pi}{M}(x+y)\right)} \sqrt{Q\left(\frac{2\pi}{M}y\right)} f(x+y)f(y) \right]^p$$

$$\leq \sum_{|x| \geq (2\pi)^{-1}Ma} \left( \sum_{y \in \mathbb{Z}^d} |f(x+y)f(y)| \right)^{p-1}$$

$$\times \sum_{y \in \mathbb{Z}^d} Q^{p/2}\left(\frac{2\pi}{M}(x+y)\right) Q^{p/2}\left(\frac{2\pi}{M}y\right) |f(x+y)f(y)|.$$

Notice that for any $x \in \mathbb{Z}^d$,

$$\sum_{y \in \mathbb{Z}^d} |f(x+y)f(y)| \leq \sum_{y \in \mathbb{Z}^d} f^2(y) = 1.$$

Thus

$$\sum_{|x| \geq (2\pi)^{-1}Ma} \left[ \sum_{y \in \mathbb{Z}^d} \sqrt{Q\left(\frac{2\pi}{M}(x+y)\right)} \sqrt{Q\left(\frac{2\pi}{M}y\right)} f(x+y)f(y) \right]^p$$

$$\leq \sum_{|x| \geq (2\pi)^{-1}Ma} \sum_{y \in \mathbb{Z}^d} Q^{p/2}\left(\frac{2\pi}{M}(x+y)\right)$$

$$\times Q^{p/2}\left(\frac{2\pi}{M}y\right) |f(x+y)f(y)|$$

$$\leq \left( \sum_{|x| \geq (2\pi)^{-1}Ma} \sum_{y \in \mathbb{Z}^d} Q^p\left(\frac{2\pi}{M}(x+y)\right) Q^p\left(\frac{2\pi}{M}y\right) \right)^{1/2}$$

$$\times \left( \sum_{x,y \in \mathbb{Z}^d} f^2(x+y)f^2(y) \right)^{1/2}$$

$$= \left( \sum_{|x-y| \geq (2\pi)^{-1}Ma} Q^p\left(\frac{2\pi}{M}x\right) Q^p\left(\frac{2\pi}{M}y\right) \right)^{1/2}.$$

Notice that

$$\frac{1}{M^{2d}} \sum_{|x-y| \geq (2\pi)^{-1}Ma} Q^p\left(\frac{2\pi}{M}x\right) Q^p\left(\frac{2\pi}{M}y\right)$$

$$\longrightarrow (2\pi)^{2d} \int\int_{\{|\lambda - \gamma| \geq a\}} Q^p(\lambda) Q^p(\gamma) \, d\lambda \, d\gamma$$

as $M \to \infty$.

For any given $\delta > 0$, therefore, one can find $a > 0$ such that

$$
\text{(A.2)} \quad \frac{1}{M^d} \sup_{|f|_2 = 1} \sum_{|x| \geq (2\pi)^{-1}Ma} \left[ \sum_{y \in \mathbb{Z}^d} \sqrt{Q\left(\frac{2\pi}{M}(x+y)\right)} \right.
$$
$$
\left. \times \sqrt{Q\left(\frac{2\pi}{M}y\right)} f(x+y)f(y) \right]^p \leq \delta
$$

for sufficiently large $M$.

For any $x = (x_1, \ldots, x_d) \in \mathbb{R}^d$, we write $[x] = ([x_1], \ldots, [x_d])$ for the lattice part of $x$. (We also use the notation $[\cdots]$ for parenthesis without causing any confusion.) For any $f \in \mathcal{L}^2(\mathbb{Z}^d)$ with $|f|_2 = 1$,

$$
\sum_{|x| \leq (2\pi)^{-1}Ma} \left[ \sum_{y \in \mathbb{Z}^d} \sqrt{Q\left(\frac{2\pi}{M}(x+y)\right)} \sqrt{Q\left(\frac{2\pi}{M}y\right)} f(x+y)f(y) \right]^p
$$
$$
= \int_{\{|\lambda| \leq (2\pi)^{-1}Ma\}} d\lambda
$$
$$
\times \left[ \int_{\mathbb{R}^d} \sqrt{Q\left(\frac{2\pi}{M}([\lambda]+[\gamma])\right)} \sqrt{Q\left(\frac{2\pi}{M}[\gamma]\right)} f([\lambda]+[\gamma])f([\gamma])\, d\gamma \right]^p
$$
$$
= \left(\frac{M}{2\pi}\right)^d \int_{\{|\lambda| \leq a\}} d\lambda
$$
$$
\times \left[ \left(\frac{M}{2\pi}\right)^d \int_{\mathbb{R}^d} \sqrt{Q_M\left(\gamma + \frac{2\pi}{M}\left[\frac{M}{2\pi}\lambda\right]\right)} \sqrt{Q_M(\gamma)} \right.
$$
$$
\left. \times f\left(\left[\frac{M}{2\pi}\lambda\right] + \left[\frac{M}{\pi}\gamma\right]\right) f\left(\left[\frac{M}{2\pi}\gamma\right]\right) d\gamma \right]^p,
$$

where

$$
Q_M(\lambda) = Q\left(\frac{2\pi}{M}\left[\frac{M}{\pi}\lambda\right]\right), \qquad \lambda \in \mathbb{R}^d.
$$

Write

$$
g_0(\lambda) = \left(\frac{M}{2\pi}\right)^{d/2} f\left(\left[\frac{M}{2\pi}\lambda\right]\right), \qquad \lambda \in \mathbb{R}^d.
$$

We have

$$
\int_{\mathbb{R}^d} g_0^2(\lambda)\, d\lambda = \left(\frac{M}{2\pi}\right)^d \int_{\mathbb{R}^d} f^2\left(\left[\frac{M}{2\pi}\lambda\right]\right) d\lambda
$$
$$
= \int_{\mathbb{R}^d} f^2([\lambda])\, d\lambda = \sum_{x \in \mathbb{Z}^d} f^2(x) = 1.
$$

We can also see that under this correspondence,

$$\left(\frac{M}{2\pi}\right)^{d/2} f\left(\left[\frac{M}{2\pi}\lambda\right] + \left[\frac{M}{2\pi}\gamma\right]\right) = g_0\left(\gamma + \frac{2\pi}{M}\left[\frac{M}{2\pi}\lambda\right]\right), \qquad \lambda, \gamma \in \mathbb{R}^d.$$

In view of (A.2), therefore, we need only to show that for any fixed $a > 0$

$$
\begin{aligned}
\text{(A.3)} \quad \limsup_{M\to\infty} \sup_{\|g\|_2=1} &\int_{\{|\lambda|\le a\}} d\lambda \left[\int_{\mathbb{R}^d} \sqrt{Q_M\left(\gamma + \frac{2\pi}{M}\left[\frac{M}{2\pi}\lambda\right]\right)} \sqrt{Q_M(\gamma)}\right. \\
&\qquad\qquad\qquad\qquad \left. \times g\left(\gamma + \frac{2\pi}{M}\left[\frac{M}{2\pi}\lambda\right]\right) g(\gamma)\, d\gamma\right]^p \\
&\le \sup_{\|g\|_2=1} \int_{\{|\lambda|\le a\}} d\lambda \left[\int_{\mathbb{R}^d} \sqrt{Q(\lambda+\gamma)} \sqrt{Q(\gamma)} g(\lambda+\gamma) g(\gamma)\, d\gamma\right]^p.
\end{aligned}
$$

Indeed, by the inverse Fourier transformation the function

$$U_M(\lambda) = \int_{\mathbb{R}^d} \sqrt{Q_M(\gamma+\lambda)} \sqrt{Q_M(\gamma)} g(\gamma+\lambda) g(\gamma)\, d\gamma$$

is the Fourier transform of the function

$$
\begin{aligned}
\text{(A.4)} \quad V_M(x) &= \frac{1}{(2\pi)^d} \int_{\mathbb{R}^d} U_M(\lambda) e^{-i\lambda\cdot x}\, d\lambda \\
&= \frac{1}{(2\pi)^d} \int_{\mathbb{R}^d} e^{-i\lambda\cdot x}\, d\lambda \int_{\mathbb{R}^d} \sqrt{Q_M(\gamma+\lambda)} \sqrt{Q_M(\gamma)} g(\gamma+\lambda) g(\gamma)\, d\gamma \\
&= \frac{1}{(2\pi)^d} \int\int_{\mathbb{R}^d\times\mathbb{R}^d} e^{-i(\lambda-\gamma)\cdot x} \sqrt{Q(\lambda)} g(\lambda) \sqrt{Q(\gamma)} g(\gamma)\, d\lambda\, d\gamma \\
&= \frac{1}{(2\pi)^d} \left|\int_{\mathbb{R}^d} e^{ix\cdot\gamma} \sqrt{Q_M(\gamma)} g(\gamma)\, d\gamma\right|^2.
\end{aligned}
$$

Therefore

$$
\begin{aligned}
\text{(A.5)} \quad &\int_{\mathbb{R}^d} \sqrt{Q_M\left(\gamma + \frac{2\pi}{M}\left[\frac{M}{2\pi}\lambda\right]\right)} \sqrt{Q_M(\gamma)}\, g\left(\gamma + \frac{2\pi}{M}\left[\frac{M}{2\pi}\lambda\right]\right) g(\gamma)\, d\gamma \\
&= U_M\left(\frac{2\pi}{M}\left[\frac{M}{2\pi}\lambda\right]\right) \\
&= \frac{1}{(2\pi)^d} \int_{\mathbb{R}^d} \exp\left\{ix \cdot \frac{2\pi}{M}\left[\frac{M}{2\pi}\lambda\right]\right\} \left|\int_{\mathbb{R}^d} e^{ix\cdot\gamma} \sqrt{Q_M(\gamma)} g(\gamma) d\gamma\right|^2 dx \\
&\le \frac{1}{(2\pi)^d} \int_{\mathbb{R}^d} \left|1 - \exp\left\{ix \cdot \left(\lambda - \frac{2\pi}{M}\left[\frac{M}{2\pi}\lambda\right]\right)\right\}\right| \\
&\qquad\qquad \times \left|\int_{\mathbb{R}^d} e^{ix\cdot\gamma} \sqrt{Q_M(\gamma)} g(\gamma)\, d\gamma\right|^2 dx
\end{aligned}
$$

$$+ \frac{1}{(2\pi)^d} \int_{\mathbb{R}^d} e^{ix\cdot\lambda} \left| \int_{\mathbb{R}^d} e^{ix\cdot\gamma} \sqrt{Q_M(\gamma)} g(\gamma) \, d\gamma \right|^2 dx.$$

By Parseval's identity and by the fact that $Q_M \leq 1$,

$$\frac{1}{(2\pi)^d} \int_{\mathbb{R}^d} \left| \int_{\mathbb{R}^d} e^{ix\cdot\gamma} \sqrt{Q_M(\gamma)} g(\gamma) \, d\lambda \right|^2 dx$$
$$= \int_{\mathbb{R}^d} Q_M(\gamma) g^2(\gamma) \, d\gamma \leq \int_{\mathbb{R}^d} g^2(\gamma) \, d\gamma = 1.$$

Hence, the first term on the right-hand side of (A.5) tends to 0 uniformly over $\lambda \in \mathbb{R}^d$ and over all $g \in \mathcal{L}^2(\mathbb{R}^d)$ with $\|g\|_2 = 1$ as $M \to \infty$. The second term on the right-hand side of (A.5) is equal to

$$\int_{\mathbb{R}^d} e^{ix\cdot\lambda} V_M(x) \, dx = U_M(\lambda) = \int_{\mathbb{R}^d} \sqrt{Q_M(\lambda+\gamma)} \sqrt{Q_M(\gamma)} g(\lambda+\gamma) g(\gamma) \, d\gamma.$$

Consequently, we will have (A.3) if we can prove

$$\begin{aligned}
(A.6) \quad & \limsup_{M\to\infty} \sup_{\|g\|_2=1} \int_{\{|\lambda|\leq a\}} d\lambda \left[ \int_{\mathbb{R}^d} \sqrt{Q_M(\lambda+\gamma)} \sqrt{Q_M(\gamma)} g(\lambda+\gamma) g(\gamma) \, d\gamma \right]^p \\
& \leq \sup_{\|g\|_2=1} \int_{\{|\lambda|\leq a\}} d\lambda \left[ \int_{\mathbb{R}^d} \sqrt{Q(\lambda+\gamma)} \sqrt{Q(\gamma)} g(\lambda+\gamma) g(\gamma) \, d\gamma \right]^p.
\end{aligned}$$

By uniform continuity of the function $Q$ we have that $Q_M(\cdot) \to Q(\cdot)$ uniformly on $\mathbb{R}^d$. Thus, given $\varepsilon > 0$ we have

$$\sup_{\lambda,\gamma\in\mathbb{R}^d} |\sqrt{Q_M(\lambda+\gamma)} \sqrt{Q_M(\gamma)} - \sqrt{Q(\lambda+\gamma)} \sqrt{Q(\gamma)}| < \varepsilon$$

for sufficiently large $M$. Therefore,

$$\begin{aligned}
& \left\{ \int_{\{|\lambda|\leq a\}} d\lambda \left[ \int_{\mathbb{R}^d} \sqrt{Q_M(\lambda+\gamma)} \sqrt{Q_M(\gamma)} g(\lambda+\gamma) g(\gamma) \, d\gamma \right]^p \right\}^{1/p} \\
& \leq \varepsilon \left\{ \int_{\{|\lambda|\leq a\}} d\lambda \left[ \int_{\mathbb{R}^d} g(\lambda+\gamma) g(\gamma) \, d\gamma \right]^p \right\}^{1/p} \\
& + \left\{ \int_{\{|\lambda|\leq a\}} d\lambda \left[ \int_{\mathbb{R}^d} \sqrt{Q(\lambda+\gamma)} \sqrt{Q(\gamma)} g(\lambda+\gamma) g(\gamma) \, d\gamma \right]^p \right\}^{1/p}.
\end{aligned}$$

Finally, (A.6) follows from the fact that

$$\int_{\{|\lambda|\leq a\}} d\lambda \left[ \int_{\mathbb{R}^d} g(\lambda+\gamma) g(\gamma) \, d\gamma \right]^p \leq C_d a^d,$$

where $C_d$ is the volume of a $d$-dimensional unit ball. $\quad\square$

Lemma A.2.    *Assume $d < 2\alpha$. Let $M_\psi$ be defined in* (1.10) *and let $\rho$ be defined in* (1.4) *with $p = 2$. Then*

$$(A.7) \qquad\qquad M_\psi = (2\pi)^{-d\alpha/(2\alpha-d)} \rho^{\alpha/(2\alpha-d)}.$$

Proof.    Replace $f(\lambda)$ by $\sqrt{Q(\lambda)} f(\lambda)$ in (1.4). Then

$$\rho = \sup_{\|f\|_{\mathcal{L}^2(Q)}=1} \int_{\mathbb{R}^d} \left[ \int_{\mathbb{R}^d} Q(\lambda+\gamma) f(\lambda+\gamma) Q(\gamma) f(\gamma) \, d\gamma \right]^2 d\lambda,$$

where

$$\|f\|_{\mathcal{L}^2(Q)} = \left( \int_{\mathbb{R}^d} f^2(\lambda) Q(\lambda) \, d\lambda \right)^{1/2}.$$

By the inverse Fourier transformation and by a computation similar to the one given in (A.4), the function

$$U(\lambda) = \int_{\mathbb{R}^d} Q(\lambda+\gamma) f(\lambda+\gamma) Q(\gamma) f(\gamma) \, d\gamma$$

is the Fourier transform of the function

$$V(x) = \frac{1}{(2\pi)^d} \left| \int_{\mathbb{R}^d} e^{-ix\cdot\gamma} Q(\gamma) f(\gamma) \, d\gamma \right|^2.$$

By Parseval's identity

$$\int_{\mathbb{R}^d} \left[ \int_{\mathbb{R}^d} Q(\lambda+\gamma) f(\lambda+\gamma) Q(\gamma) f(\gamma) \, d\gamma \right]^2 d\lambda$$

$$= (2\pi)^d \int_{\mathbb{R}^d} V^2(x) \, dx$$

$$= \frac{1}{(2\pi)^d} \int_{\mathbb{R}^d} \left| \int_{\mathbb{R}^d} e^{ix\cdot\gamma} Q(\gamma) f(\gamma) \, d\gamma \right|^4 dx.$$

Let $p_t(x)$ be the density of $X(t)$ and write

$$G(x) = \int_0^\infty p_t(x) e^{-t} \, dt, \qquad x \in \mathbb{R}^d.$$

Notice that

$$\int_{\mathbb{R}^d} e^{i\lambda\cdot x} G(x) \, dx = Q(\lambda).$$

If we consider $f(\lambda)$ as the Fourier transform of the function $g(x)$ on $\mathbb{R}^d$, then

$$\int_{\mathbb{R}^d} e^{-ix\cdot\gamma} Q(\gamma) f(\gamma) \, d\gamma = (2\pi)^d \int_{\mathbb{R}^d} G(y-x) g(y) \, dy = (2\pi)^d Gg(x),$$

$$\|f\|_{\mathcal{L}^2(Q)}^2 = (2\pi)^d \int_{\mathbb{R}^d \times \mathbb{R}^d} G(y-x) g(x) g(y) \, dx \, dy = (2\pi)^d \langle g, Gg \rangle.$$

Summarizing the above steps, we obtain

$$(A.8) \qquad \rho = (2\pi)^{3d} \sup_{\langle g, Gg \rangle = (2\pi)^{-d}} \int_{\mathbb{R}^d} |Gg(x)|^4 \, dx.$$

Write $h(x) = Gg(x)$ and recall the resolvent identity

$$I = G - \mathcal{A} \circ G,$$

where $I$ is identity operator and where $\mathcal{A}$ is the infinitesimal generator of the Markov process $X(t)$. Then

$$\langle g, Gg \rangle = \langle h - \mathcal{A}h, h \rangle = \|h\|_2 + \int_{\mathbb{R}^d} \psi(\lambda) |\widehat{h}(\lambda)|^2 \, d\lambda = \|h\|_2 + \|\widehat{h}\|^2_{\mathcal{L}^2(\psi)},$$

where

$$\|f\|^2_{\mathcal{L}^2(\psi)} = \int_{\mathbb{R}^d} \psi(\lambda) |f(\lambda)|^2 \, d\lambda$$

and where the second step follows from the fact (page 24 in [3]) that

$$\langle \mathcal{A}h, h \rangle = -\int_{\mathbb{R}^d} \psi(\lambda) |\widehat{h}(\lambda)|^2 \, d\lambda.$$

Hence, from (A.8) we have

$$(A.9) \qquad \rho = (2\pi)^{3d} \sup_{\|h\|_2 + \|\widehat{h}\|^2_{\mathcal{L}^2(\psi)} = (2\pi)^{-d}} \int_{\mathbb{R}^d} |h(x)|^4 \, dx.$$

Write

$$M_\psi(\theta) = \sup_{g \in \mathcal{F}_\psi} \left\{ \theta \left( \int_{\mathbb{R}^d} |g(x)|^4 \, dx \right)^{1/2} - \int_{\mathbb{R}^d} \psi(\lambda) |\widehat{g}(\lambda)|^2 \, d\lambda \right\}, \qquad \theta > 0,$$

where $\mathcal{F}_\psi$ is defined in (1.11). By (2.10) in [8] (with $p = 2$),

$$(A.10) \qquad M_\psi(\theta) = \theta^{2\alpha/(2\alpha-d)} M_\psi, \qquad \theta > 0.$$

Therefore, we will have (A.7) if we can prove that

$$(A.11) \qquad M_\psi \left( \frac{(2\pi)^{d/2}}{\sqrt{\rho}} \right) = 1.$$

Indeed, for any $0 < \varepsilon < \rho$, by (A.9) there is an $h_0$ such that $\|h_0\|_2 + \|\widehat{h_0}\|^2_{\mathcal{L}^2(\psi)} = (2\pi)^{-d}$ and that

$$\int_{\mathbb{R}^d} |h_0(x)|^4 \, dx > (2\pi)^{-3d}(\rho - \varepsilon).$$

Consequently,

$$M_\psi\left(\frac{(2\pi)^{d/2}}{\sqrt{\rho}}\right) \geq \frac{(2\pi)^{d/2}/\sqrt{\rho-\varepsilon}(\int_{\mathbb{R}^d}|h_0(x)|^4\,dx)^{1/2} - \int_{\mathbb{R}^d}\psi(\lambda)|\widehat{h_0}(\lambda)|^2\,d\lambda}{\int_{\mathbb{R}^d}|h_0(x)|^2\,dx}$$

$$\geq \frac{(2\pi)^{-d} - \int_{\mathbb{R}^d}\psi(\lambda)|\widehat{h_0}(\lambda)|^2\,d\lambda}{\int_{\mathbb{R}^d}|h_0(x)|^2\,dx} = 1.$$

Let $\varepsilon \to 0^+$ on the left-hand side. By (A.10), $M(\theta)$ is continuous. So we have

$$(A.12) \qquad\qquad M_\psi\left(\frac{(2\pi)^{d/2}}{\sqrt{\rho}}\right) \geq 1.$$

On the other hand, by (A.9) again

$$M_\psi\left(\frac{(2\pi)^{d/2}}{\sqrt{\rho}}\right)$$

$$= \sup_{g\in\mathcal{F}_\psi}\left\{\frac{(2\pi)^{d/2}}{\sqrt{\rho}}\left(\int_{\mathbb{R}^d}|g(x)|^4\,dx\right)^{1/2} - \int_{\mathbb{R}^d}\psi(\lambda)|\widehat{g}(\lambda)|^2\,d\lambda\right\}$$

$$\leq \sup_{g\in\mathcal{F}_\psi}\left\{\frac{(2\pi)^{d/2}}{\sqrt{\rho}}(2\pi)^{-d/2}\sqrt{\rho}\left[1 + \int_{\mathbb{R}^d}\psi(\lambda)|\widehat{g}(\lambda)|^2\,d\lambda\right]\right.$$

$$\left. - \int_{\mathbb{R}^d}\psi(\lambda)|\widehat{g}(\lambda)|^2\,d\lambda\right\} = 1. \qquad \square$$

**Acknowledgments.** I thank Peter Mörters for helpful discussions on the combinatorial method of high moment asymptotics, and an anonymous referee for his/her comments and suggestions.

DEPARTMENT OF MATHEMATICS
UNIVERSITY OF TENNESSEE
KNOXVILLE, TENNESSEE 37996
USA
E-MAIL: xchen@math.utk.edu
URL: http://www.math.utk.edu/~xchen